\documentclass[12pt,leqno]{article}
\usepackage[francais]{babel}
\usepackage[latin1]{inputenc}
\usepackage[T1]{fontenc}
\usepackage{amsmath,amsfonts,amssymb,amsthm,amscd}
\usepackage{bm}
\date{}
\newcommand\1{\mbox{1\hspace{-.3em}I}}
\newcommand{\ba}{\backslash}
\newcommand\rg{\rightarrow}
\newcommand\lgr{\longrightarrow}  

\newcommand\N{\mathbb{N}}

\newcommand\R{\mathbb{R}}

\newcommand\CA{\mathcal{A}}
\newcommand\Cc{\mathcal{C}}

\newcommand\Nn{\mathcal{N}}
\newcommand\CPp{\mathcal{P}}
\newcommand\Dd{\mathcal{D}}
\newcommand\Ee{\mathcal{E}}

\newcommand\Cs{\mathcal{S}}

\newcommand\inv{\mathrm{inv}}

\newcommand{\newsmile}{\smallsmile_{\hspace{-3mm}\displaystyle\smallfrown}}

\def\adots{\mathinner{\mkern1mu\raise1pt\vbox{\kern7pt\hbox{.}}
\mkern2mu\raise4pt\hbox{.}
\mkern2mu\raise7pt\hbox{.}\mkern1mu}} 

\newtheorem{montheo}{\textsc{Th\'eor\`eme}}[section]
\newtheorem{monlem}{\textsc{Lemme}}[section]

\newtheorem{maprop}[monlem]{\textsc{Proposition}}

\begin{document}

\title{Etude du graphe divisoriel 5}

\author{Eric Saias}

%\date{juillet 2015}

\maketitle

\begin{center}
Pour Michel Balazard,

à l'occasion de son soixante et unième anniversaire
\end{center}

\section{Introduction}

On appelle graphe divisoriel le graphe non orienté sur $\N^*$ où deux entiers sont en  relation quand le petit divise le grand. On appelle chaîne de longueur $\ell$ tout chemin simple du graphe divisoriel dont le nombre de sommets est $\ell$. On note $f(x)$ la longueur maximum d'une chaîne de la restriction du graphe divisoriel aux entiers $\le x$. 

C'est en lien avec le travail  \cite{ErFrHe} d'Erd\"os, Freud et Hegyv\'ari, que Pomerance \cite{Pom} introduit la fonction $f(x)$ et résout une conjecture de Hegyv\'ari en montrant que $f(x)=o(x)$. Avec la minoration de Pollington \cite{Pol}, on sait en 1983 qui est l'année de parution des trois articles \cite{ErFrHe}, \cite{Pom} et \cite{Pol} que
$$
\frac{x}{\exp\{(2+o(1))\sqrt{\log x\, \log\log x}\}} \le f(x) = o(x)\,, \quad (x\rg +\infty).
$$
Les travaux de Tenenbaum \cite{Ten2} et moi--même \cite{Sai1}, \cite{Sai3} et \cite{Sai4} ont permis d'établir dix ans plus tard qu'il existe deux réels $a$ et $b$ avec $b>a>0$ tels que pour tout $x\ge2$
\begin{equation}
a\ \frac{x}{\log x} \le
f(x) \le b\ \frac{x}{\log x}\,.
\end{equation}
Aucune valeur explicite de $a$ et $b$ n'avait été calculée jusqu'à présent, même pour les $x$ suffisamment grands.
Des considérations théoriques et des calculs numériques m'amènent à penser que l'on a 

\eject

\noindent\textsc{Conjecture F}.

\textit{Il existe un réel $c$ vérifiant $3<c<7$ tel que
}
$$
f(x) \sim c\ \frac{x}{\log x}, \qquad (x\rg +\infty)\,.
$$

\vspace{2mm}

Découpons la recherche de cette conjecture en trois questions distinctes : montrer que $f(x)/(x/\log x)$ converge quand $x\rg +\infty$,

\vspace{2mm}
\begin{tabular}{l l}
\hspace{-9mm} que &$\alpha:= \liminf\limits_{x\rg +\infty} (f(x)/(x/\log x)) >3$,\\
\vspace{4mm}
\hspace{-7mm}et que &$\beta:= \limsup\limits_{x\rg +\infty}(f(x)/(x/\log x)) <7$.
\end{tabular}

\vspace{2mm}

Nous allons nous focaliser ici sur la minoration de $\alpha$, que l'on peut aborder maintenant grâce au récent travail numérique de Weingartner sur les entiers à diviseurs denses \cite{Wei5}.

Une question plus profonde est celle de la structure des chaînes de longueur maximum. Elle est conjecturalement approximativement fractale. Voyons cela.

Chadozeau \cite{Cha} a montré que $f(100)=77$. Donnons un exemple de chaîne d'entiers $\le 100$ formée de $77$ entiers.

\vspace{2mm}

\noindent $62-31-93-1-87-29-58-2-92-46-23-69-3-$

\noindent $57-19-38-76-4-68-34-17-85-5-65-13-52-26-78-$

\noindent $6-24-48-96-12-72-36-18-54-27-81-9-$

\noindent $63-21-42-84-28-56-14-98-49-7-35-70-$

\noindent $10-60-30-90-45-15-75-25-50-100-20-40-80-$

\noindent $16-32-64-8-88-44-22-66-33-99-11-55$.

\vspace{2mm}

On constate que les entiers ayant le même plus grand facteur premier ont tendance à être regroupés en sous--chaînes. Par exemple ceux dont le plus grand facteur premier est $7$ sont tous les entiers de la quatrième ligne. Pour faire apparaître cette structure, on récrit cette chaîne dans la forme suivante où on factorise une sous--chaîne à plus grand facteur premier constant par ce plus grand facteur premier. 

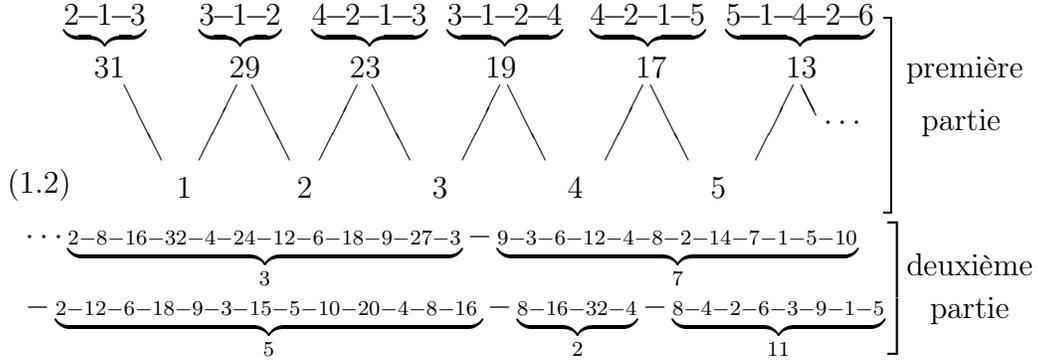
\eject

%\vspace{2mm}
\begin{equation*}
\unitlength=1cm
\begin{picture}(5,2)
\put(-3.6,1.6){$\underbrace{2\!\!-\!\!1\!\!-\!\!3}$}
\put(-3.2,0.9){$31$}
\put(-2.8,0.8){\line(1,-2){0.5}}
\put(-1.8,1.6){$\underbrace{3\!\!-\!\!1\!\!-\!\!2}$}
\put(-1.4,0.9){$29$}
\put(-1.8,-0.2){\line(1,2){0.5}}
\put(-2.1,-0.7){$1$}
\put(-1.1,0.8){\line(1,-2){0.5}}
\put(-0.5,-0.7){$2$}
\put(-0.3,1.6){$\underbrace{4\!\!-\!\!2\!\!-\!\!1\!\!-\!\!3}$}
\put(0.2,0.9){$23$}
\put(-0.2,-0.2){\line(1,2){0.5}}
\put(0.5,0.8){\line(1,-2){0.5}}
\put(1.3,-0.7){$3$}
\put(1.5,1.6){$\underbrace{3\!\!-\!\!1\!\!-\!\!2\!\!-\!\!4}$}
\put(2,0.9){$19$}
\put(1.7,-0.2){\line(1,2){0.5}}
\put(2.3,0.8){\line(1,-2){0.5}}
\put(3.1,-0.7){$4$}
\put(3.4,1.6){$\underbrace{4\!\!-\!\!2\!\!-\!\!1\!\!-\!\!5}$}
\put(4,0.9){$17$}
\put(3.6,-0.2){\line(1,2){0.5}}
\put(4.2,0.8){\line(1,-2){0.5}}
\put(5,-0.7){$5$}
\put(5.2,1.6){$\underbrace{5\!\!-\!\!1\!\!-\!\!4\!\!-\!\!2\!\!-\!\!6}$}
\put(6,0.9){$13$}
\put(5.6,-0.2){\line(1,2){0.5}}
\put(6.2,0.8){\line(1,-2){0.2}}
\put(6.5,0.2){$\cdots$}
\put(7.6,0.9){première}
\put(7.8,0.2){partie}
\put(7.4,-0.9){\line(0,1){2.6}}
\put(7.4,1.7){\line(-1,0){0.1}}
\put(7.4,-0.9){\line(-1,0){0.1}}
\end{picture}
\end{equation*}
\noindent(1.2)
$$
\left.
\begin{array}{l}
\cdots \underbrace{\scriptstyle 2-8-16-32-4-24-12-6-18-9-27-3}_3-\underbrace{\scriptstyle 9-3-6-12-4-8-2-14-7-1-5-10}_7 \\
-\underbrace{\scriptstyle 2-12-6-18-9-3-15-5-10-20-4-8-16\,}_5 
-\underbrace{\scriptstyle 8-16-32-4}_2 -\underbrace{\scriptstyle 8-4-2-6-3-9-1-5}_{11}
\end{array}
\!\!\right]\!\!{\textrm{deuxième}\atop \textrm{partie}}
$$

Remarquons que dans la première partie de cette chaîne qui va de 62 à 78, les sous--chaînes à plus grand facteur premier constant $p$ (étage du haut) sont séparées par les \og connecteurs extérieurs\fg\ que sont les entiers 1, 2, 3, 4 et 5 (étage du bas). Ces connecteurs sont nécessaires car ces nombres premiers $p$ sont $>10=\sqrt{100}$. 

L'étude de la deuxième partie nous amène à appeler connecteur intérieur un entier à l'extrémité d'une sous--chaîne $p\Cc(x/p,p)$, qui divise l'extrémité d'une sous--chaîne $\Cc(x/p',p')$ voisine avec $p<p'$. Ces connecteurs intérieurs sont $6$, $9$, $10$, $16$ et $8$. 

Au theorem 1 de \cite{ErFrHe}, Erd\"os, Freud et Hegyv\'ari ont introduit la fonction $L(N)$ qui représente le minimum sur toutes les permutations $\sigma$ de $\{1,2,\ldots,N\}$ de la quantité $\max\limits_{1\le k<N} ppcm(\sigma(k),\sigma(k+1))$ et  ont montré que $L(N) \sim N^2/4\log N$. Pour la majoration, ils ont été amené à construire une permutation $\sigma$ pour laquelle le $N$--uplet $(\sigma(1),\sigma(2),\ldots,\sigma(N))$ a une structure fortement apparentée à celle de la chaîne (1.2). En particulier, ils regroupent les entiers ayant le même plus grand facteur premier $p$ en des intervalles qu'ils appellent blocs, et ils séparent ces blocs par des connecteurs extérieurs quand $p>\sqrt{N}$. L'analogie est naturelle car les deux problèmes sont apparentés. C'est d'ailleurs le travail \cite{ErFrHe} qui est à l'origine de tous les travaux sur le graphe divisoriel.

De manière générale, on désigne par $P(n)$ le plus grand facteur premier de l'entier  $n\ge2$ et on note par convention $P(1)=1$. On note alors $\Cs(x,y)=\{n\le x : P(n)\le y\}$ et $\Psi(x,y)=|\Cs(x,y)|$. On utilisera toujours également la lettre $p$ pour désigner un nombre premier générique.

Notons $\Cc(x,y)$ une chaîne générique d'entiers de $\Cs(x,y)$. On étend $f(x)$ en la fonction $f(x,y)$ de deux variables qui représente la longueur maximum d'une telle chaîne~$\Cc(x,y)$.

Au vu de l'exemple (1.2) de chaîne de longueur maximum $f(100)$, une stratégie naturelle pour minorer $f(x)$ est de minorer de manière itérative $f(x,y)$ par la construction d'une chaîne $\Cc(x,y)$ de la forme

\setcounter{equation}{2}
\begin{equation}
\Cc(x,y) :
\Bigg|
\begin{array}{ccl}
p_1\Cc (\frac{x}{p_1},p_1) & -\cdots - &p_j \Cc(\frac{x}{p_j},p_j)-\\
p_{j+1}\Cc(\frac{x}{p_{j+1}},p_{j+1})& - \cdots - & p_k \Cc(\frac{x}{p_k},p_k)
\end{array}
\end{equation}
où $p_1,p_2,\ldots,p_k$ sont $k$ nombres premiers $\le y$ deux à deux distincts. La plupart des chaînes qui ont été construites à partir de 1992 pour minorer $f(x)$ et ses variantes ont effectivement la structure (1.3) (voir \cite{Sai1}, \cite{Ten2}, \cite{Sai2}, \cite{Sai6} et \cite{MeSa}), quitte à introduire des connecteurs extérieurs, comme dans la première partie de (1.2). C'est encore le cas pour les deux chaînes que l'on construit dans ce présent travail.

C'est toujours en lien avec l'exemple (1.2) qu'il est naturel de conjecturer que l'on a
\begin{equation}
f(x,y) = (1+o(1)) \sum_{p\le \min(y,\sqrt{x})}
f(x/p,p)
\end{equation}
quand $x$ tend vers $+\infty$, uniformément en $y$. Cette identité de Buchstab approchée conjecturale de $f(x,y)$ est l'une des manières de comprendre le lien entre les longues chaînes et les entiers à diviseurs denses, par l'intermédiaire de la fonction de Schinzel--Szekeres. Voyons cela.

Pour tout entier $n\ge 2$, on désigne par  $1=d_1(n)<d_2(n)<\cdots<d_{\tau(n)}(n)=n$ la suite croissante de ses diviseurs. On dit que l'entier $n$ est à diviseurs $t$--denses quand il vérifie
$$
\max_{1\le i<\tau(n)} \frac{d_{i+1}(n)}{d_i(n)} \le t\,.
$$
C'est Erdös \cite{Erd} qui s'est le premier intéressé à ce type d'entiers en 1948. Ils ont été ensuite étudiés notamment par Tenenbaum et moi--même, entre 1986 et 1999 (\cite{Ten1}, \cite{Ten2}, \cite{Sai3}, \cite{Sai5}). Les remarquables travaux de Weingartner depuis 2001, dont un avec Pomerance, et un avec Pomerance et Thompson, ont énormément amélioré la connaissance que nous avons de ces entiers  (\cite{Wei1}, \cite{Wei2}, \cite{Wei3} \cite{Wei4}, \cite{PoThWe}, \cite{Wei5}, \cite{PoWe}, \cite{Wei6}, \cite{Wei7} et \cite{Wei8}).

Par convention, on écrira que pour $n=1$
$$
\max_{1\le i<\tau(n)} \frac{d_{i+1}(n)}{d_i(n)}=1.
$$

Notons
$$
\begin{array}{ll}
D(x,t)& = \Big|\Big\{n\le x : \max\limits_{1\le i <\tau(n)} \frac{d_{i+1}(n)}{d_i(n)} \le t\Big\}\Big|\\
\mathrm{et}&\\
D'(x,t)& = \Big|\Big\{n\le x : \max\limits_{1\le i <\tau(n)} \frac{d_{i+1}(n)}{d_i(n)} \le t \ \mathrm{et}\ \mu^2(n)=1\Big\}\Big|
\end{array}
$$
où  $\mu(n)$ désigne la fonction de Möbius. En 2014, Weingartner (Corollary 1.1 de \cite{Wei4}) a montré que pour tout réel $t\ge 2$ fixé, il existe un nombre réel $c_t>0$ tel que
\begin{equation}
D(x,t) = c_t \frac{x}{\log x} \Big(1+O_t \Big(\frac{1}{\log x}\Big)\Big),\qquad (x\ge 2).
\end{equation}
Deux ans plus tard, Pomerance, Thompson et Weingartner ont établi un résultat analogue pour $D'(x,t)$ (cela découle du theorem~3.1 de \cite{PoThWe}) : il existe un nombre réel $c_t'>0$ tel que 
\begin{equation}
D'(x,t) = c_t' \frac{x}{\log x} \Big(1+O_t \Big(\frac{1}{\log x}\Big)\Big),\qquad (x\ge 2).
\end{equation}
De plus, Weingartner a démontré récemment un encadrement pour certains de ces réels $c_t$ et $c_t'$ (\cite{Wei5} Corollary~2 et Corollary~7). Il  obtient en particulier
\begin{equation}
c_2 = 1,2248\cdots\\
\end{equation}
et
\begin{equation}
c_2' = 0,0686\cdots
\end{equation}

Indépendamment, pour étudier une autre question soulevée par Erd\"os, Schinzel et Szekeres \cite{ScSz} ont introduit en 1959 la fonction 
\begin{equation*}
S(n) = \Bigg|
\begin{array}{cc}
\max\limits_{1\le k\le \Omega(n)}& p_1(n)p_2(n)\cdots p_{k-1}(n) p_k^2(n) \qquad \mathrm{si} \ n\ge 2\\
\\
1 &\mathrm{si}\ n=1\,.
\end{array}
\end{equation*}
où $n=p_1(n)p_2(n)\cdots p_{\Omega(n)}(n)$ est la décomposition de $n$ en facteurs premiers avec $p_1(n)\ge p_2(n) \ge\cdots\ge p_{\Omega(n)}(n)$. On désigne par $P^-(n)$ le plus petit facteur premier de l'entier $n\ge 2$ et on note par convention $P^-(1)=+\infty$.

Notons 
\begin{equation*}
\begin{array}{rl}
\CA(x,y)&=\{n:P(n)\le y, S(n)\le x\},\\
\CA'(x,y)&=\{n:P(n)\le y, S(n)\le x\ \mathrm{et}\ \mu^2(n)=1\}\\
\CA(x,y,z)&=\{n:P(n)\le y, P^-(n)>z,\ S(n)\le x\},
\end{array}
\end{equation*}
$\CA(x)=\CA(x,x)$ et $\CA'(x)=\CA'(x,x)$. On note également les cardinaux $A(x,y)=|\CA(x,y)|$, $A'(x,y)=|\CA'(x,y)|$, $A(x,y,z)=|\CA(x,y,z)|$, $A(x)=|\CA(x)|$ et $A'(x)=|\CA'(x)|$.

La fonction $A(x,y)$ vérifie l'identité de Buchstab
$$
A(x,y) = 1+\sum_{p\le \min(y,\sqrt{x})}A(x/p,p),\qquad (x\ge 1, y\ge2)
$$
qui correspond à la formule conjecturale (1.4) pour $f(x,y)$, une fois que l'on a supprimé le terme d'erreur $o(1)$. Cela constitue une manière de voir 
 le lien entre les longues chaînes d'entiers de $\Cs(x,y)$ avec l'ensemble $\CA(x,y)$. Par ailleurs Tenenbaum a montré (lemme~2.2 de \cite{Ten1}) que
\begin{equation}
S(n) = n \max_{1\le i <\tau(n)}\frac{d_{i+1}(n)}{d_i(n)},
\end{equation}
ce qui explicite cette fois le lien entre entiers à diviseurs denses et la fonction de Schinzel--Szekeres $S(n)$

Dans \cite{Ten2}, Tenenbaum utilise ces liens pour montrer que\break $f(x)=(x/\log x)\cdot (\log\!\log x)^{O(1)}$. Pour la minoration, il construit une chaîne $\Cc(x,y)$ d'entiers de $\Cs'(x,y):=\{n\in \Cs(x,y):\mu^2(n)=1\}$ de la forme (1.3) qui vérifie
$$
\Cc(x,y) \supset \CA'(x/2,y).
$$
On en déduit immédiatement que
\begin{equation}
f(x,y) \ge A'(x/2,y), \qquad (x\ge 1, y\ge 2).
\end{equation}

Par une variante de la construction de Tenenbaum, on peut lever la contrainte pour ces chaînes d'être constituées d'entiers sans facteur carré. Autrement dit on construit ici une chaîne $\Cc(x,y)$ ayant la forme (1.3) et qui vérifie
$$
\Cc(x,y) \supset \CA(x/2,y).
$$
D'où le

\vspace{2mm}
{\montheo On a pour tout $x\ge 1$ et $y\ge 2$
\begin{equation}
f(x,y) \ge A(x/2,y).
\end{equation}}

\vspace{2mm}

Il découle immédiatement de (1.9) que $A(x) \ge D(x/2,2)$ et $A'(x)\ge D'(x/2,2)$. Avec (1.6) et (1.8), le résultat de Tenenbaum (1.10) fournit donc la minoration
$$
\alpha>0,017.
$$
Avec (1.5), (1.7) et le théorème 1, on obtient cette fois
\begin{equation}
\alpha
>0,306.
\end{equation}

On a donc gagné un facteur dix--huit. Il reste encore à gagner un facteur dix pour établir la conjecture~$\alpha>3$.

Etudions à présent le cas où $y$ est petit devant $x$. En 2008, Mazet (théorème~3.0.2 de \cite{Maz}) a montré que pour tout $y\ge 2$ fixé, dès que $x$ est suffisamment grand relativement à $y$, la restriction du graphe divisoriel à $S(x,y)$ est un graphe Hamiltonien. Autrement dit pour tout $y$ fixé, $f(x,y)=\Psi(x,y)$ pour $x\ge x_0(y)$. Le théorème~1.1 fournit uniquement $f(x,y) \sim \Psi(x,y)$ mais dès que $y=o(\log x)$. Plus précisément on a le

\vspace{2mm}

{\montheo On a uniformément pour}
\begin{equation}
2\le y\le \log x
\end{equation}
$$
\frac{\Psi(x,y)}{1+\frac{y}{\log x}} \Big(1+O\Big(\frac{\log y}{y}+\frac{y}{\log x\cdot \log y}\Big)\Big) \le f(x,y) \le \Psi(x,y).
$$

Cet encadrement de $f(x,y)$ est de même nature que celui du theorem~1.1 du très récent travail de Mehdizadeh \cite{Meh} sur la question du problème de la table de multiplication d'Erdös pour les entiers sans grand facteur premier. D'ailleurs on reprend une partie de son argumentation.

Donnons la preuve de la minoration du théorème~1.2, qui se découpe en trois étapes.
La première est d'utiliser la minoration du théorème 1.1
$$
f(x,y) \ge A(\frac{x}{2},y).
$$
La seconde consiste à utiliser la minoration immédiate
$$
A(\frac{x}{2},y) \ge \Psi(\frac{x}{2y},y).
$$
La troisième et dernière est de procéder comme dans la preuve du Corollary~2.3 de \cite{Meh}, en montrant que l'on a
$$
\Psi(\frac{x}{2y},y) = \frac{\Psi(x,y)}{1+\frac{y}{\log x}} \Big(1+O\Big(\frac{\log y}{y}+\frac{y}{\log x\cdot \log y}\Big)\Big)
$$
dans le domaine (1.13).

Dans la suite de cet article, on utilise les notations suivantes. Soient $u$ et $v$ deux fonctions définies sur une partie $D$ de $\R^{+d}$ à valeurs dans $\R^{+*}$. On note $u(X)\ll v(X)$ (respectivement $u(X)\gg v(X))$ pour signifier qu'il existe un réel positif $K$ tel que pour tout élément $X$ de $D$, on a $u(X)\le Kg(X)$ (resp. $u(X)\ge K v(X))$. On écrit $u(X) \asymp v(X)$ quand on a simultanément $u(X) \ll v(X)$ et $u(X)\gg v(X)$. Quand la constante $K$ dépend d'un paramètre $\eta$, c'est à dire est une fonction de $\eta$, on écrira par exemple $u(X)\ll_\eta v(X)$ à la place de $u(X)\ll v (X)$.

On s'intéresse à présent à une autre question relative aux longues chaînes d'entiers.

Notons $R(x,z)$ le cardinal maximum de l'union de chaînes d'entiers $\le x$ deux à deux disjointes et toutes de longueur $z$. On démontre dans \cite{Sai7} que l'on a
\begin{equation}
R(x,z) \asymp \frac{x}{\log z}, \qquad (2\le z \le f(x)).
\end{equation}
Désignons par $f_a(x)$ la longueur maximum d'une chaîne d'entiers de $\CA(x)$. Pour montrer (1.14) on met en lumière une équation fonctionnelle approchée de $R(x,z)$, qui fait intervenir la fonction $f_a(x)$. Pour itérer cette équation fonctionnelle, on a besoin de la minoration
\begin{equation}
f_a(x) \gg x/\log x.
\end{equation}

Pour démontrer cette minoration, ma première idée était ambitieuse : il s'agissait de montrer que $f_a(x)$ prend pour tout $x$ la plus grande valeur possible, i.e. $f_a(x) =A(x)$ et on a (1.15) par l'estimation $A(x)\asymp x/\log x$ (voir le lemme~4.2 ci--dessous). Autrement dit, il s'agissait de confirmer l'

\vspace{2mm}

\textsc{Hypothèse A}. \textit{Pour tout $N\ge 1$, on peut ordonner l'ensemble $\CA(N)$ en une chaîne.}

\vspace{2mm}

Il serait à désirer, sans doute, que l'on eût une démonstration rigoureuse de cette proposition ; néanmoins j'ai laissé cette recherche de côté pour le moment après quelques rapides essais infructueux car elle paraît superflue pour le but immédiat de mon étude.

En effet notons $f^*(x)$ (respectivement $f_a^*(x)$) la longueur maximum d'une chaîne d'entiers de $]\sqrt{x},x]$ (resp. et qui appartiennent à $\CA(x)$). J'ai constaté alors qu'au prix de complications techniques, on pouvait établir (1.15) par l'intermédiaire de $f_a^*(x)\gg x/\log x$, en réutilisant la méthode employée au chapitre~8 de \cite{Sai2} pour minorer $f^*(x)$. On a donc le

\vspace{2mm}

{\montheo Il existe deux réels $c$ et $K$ avec $K> c>0$ tels que pour $x \ge 2$, on~a}
$$
c\frac{x}{\log x} \le f_a(x) \le K\frac{x}{\log x}.
$$

La majoration découle immédiatement de (1.1). Pour la minoration on note $S^*(x,y)=\{n : \sqrt{x}<n \le x :P(n)\le y)$. Notons $f_a(x,y)$ (respectivement $f_a^*(x,y)$) la longueur maximum d'une chaîne d'entiers de $\CA(x) \cap \Cs(x,y)$ (resp. $\CA(x) \cap \Cs^*(x,y)$). On montre en fait le

\vskip2mm

{\montheo Il existe un réel $M\ge 1$ tel que pour tous $x\ge 8$ et $y\ge 2$, on a}
$$
f_a(x,y) \ge f_a^* (x,y) \ge A(\frac{x}{M}, y , 2)\,.
$$

\vskip2mm

On sait que $A(x,x,2) \asymp x/\log x$ (voir le lemme 4.2). Donc la minoration $f_a(x) \gg x/\log x$ du théorème~1.3 résulte immédiatement du théorème~1.4.

L'itération d'identités ou d'inéquations de Buchstab est un outil important de la théorie multiplicative des nombres. Suivant les cas on est amené à cribler par les petits ou les grands facteurs premiers. Notons $\pi(x,k)=|\{n\le x : \omega(n)=k\}|$. Dans \cite{Bal}, Balazard crible par les petits facteurs premiers en  prolongeant $\pi(x,k)$ par
$$
\pi (x,z,k):=|\{n\le x:P^-(n)>z\ \mathrm{et }\ \omega(n)=k\}|\,.
$$
C'est en itérant l'identité de Buchstab que cette dernière fonction vérifie, qu'il a prouvé la conjecture d'Erd\"os suivante : pour tout $x$ suffisamment grand, la suite $(\pi(x,k))_{k\ge 1}$ est unimodale.

L'itération d'identités de Buchstab a été également utilisée plusieurs fois dans le cadre de l'estimation de la fonction de comptage des entiers à diviseurs denses : dans \cite{Ten2}, \cite{Sai5} et pour le Lemme 3 de \cite{PoWe}, on itère l'identité de Buchstab obtenue en classant les entiers par leur plus petit facteur premier après avoir criblé par les petits facteurs premiers. Dans \cite{Sai2}, \cite{Sai3} et \cite{Sai4}, c'est avec les grands facteurs premiers. 
Notons $PP_1(x)=|\{p\le x:p-1$ est pratique$\}|$. Pomerance et Weingartner ont montré très récemment que $PP_1(x)\ll x/\log^2x$, ce qui est l'ordre de grandeur attendu (Theorem 1 de \cite{PoWe}). Pour cela, ils réutilisent la  méthode d'itération d'identité de Buchstab de \cite{Sai3}.

 Le travail \cite{Sai1} a permis de dégager une méthode générale pour minorer la longueur maximum d'une chaîne d'entiers vérifiant au moins la propriété d'appartenir à $\Cs(x,y)$. Il s'agit de combiner la construction d'une chaîne d'entiers ayant la structure (1.3) avec l'itération d'une inégalité de Buchstab. Cette méthode a donc été utilisée dans \cite{Sai1}, mais aussi dans \cite{Sai2}, aux lemmes 5 et 6 de \cite{Sai6} et dans \cite{MeSa}. C'est encore comme cela que l'on travaille ici pour établir le théorème~1.4.

\section{Opérations sur les chaînes}

Les constructions de longues chaînes s'obtiennent en effectuant diverses opérations sur des sous--chaînes (voir \cite{Sai1}, \cite{Ten2}, etc). 

Plutôt que de donner une définition formelle, on rappelle la nature de ces opérations sur des exemples. Soit $\Cc_1=2-6-3$. La chaîne $3-6-2$ est l'inverse de $\Cc_1$, et est notée inv $\Cc_1$. La chaîne $14-42-21$ est le produit de $\Cc_1$ par 7. On la note $7\Cc_1$. Il y a maintenant trois variantes de \og collage \fg{} de deux chaînes.

\vspace{2mm}

Si $\Cc_2=12-24-8$, 

\noindent on appellera concaténation de $\Cc_1$ et $\Cc_2$ et on notera $\Cc_1-\Cc_2$ la chaîne $2-6-3-12-24-8$.

\vspace{2mm}
Si $\Cc_2=3-12-24$, 

\noindent on appellera collage de $\Cc_1$ et $\Cc_2$ et on notera encore $\Cc_1-\Cc_2$ la chaîne $2-6-3-12-24$.

\vspace{2mm}

Si  $\Cc_2=5-40-20$,

\noindent on appellera juxtaposition de $\Cc_1$ et $\Cc_2$ par le connecteur extérieur $15$, la chaîne $2-6-3-15-5-40-20$. On la notera
$$
\begin{array}{ccccc}
\Cc_1&&&&\Cc_2\\
&\diagdown &&\diagup\\
&&15
\end{array}
$$
On conviendra enfin que $\emptyset$ est une chaîne.

\section{Preuve du théorème 1.1}

Notons $p_1<p_2<\cdots$ la suite croissante des nombres premiers. On pose par convention $p_0=2$ et $p_{-1}=1/2$. On note
$$\begin{array}{ll}
j=j(x,y)=\pi(\min(y,x/2))\\
\hspace{-3.8cm}\mathrm{et}&\\
k=k(x,y) =\pi (\min(y,\sqrt{x/2})).
\end{array}
$$

On démontre le théorème 1.1 sous la forme plus précise suivante :

\vspace{2mm}

{\maprop Soient $x\ge 2$ et $y\ge 2$.

\noindent Il existe une chaîne $\Cc(x,y)$ d'entiers de $\Cs(x,y)$ qui commence en $2p_{j(x,y)-1}$, s'achève en $2$ et qui contient $\CA(x/2,y)$.}

\vspace{2mm}

\noindent\textbf{Démonstration.} On construit la chaîne $\Cc(x,y)$ par récurrence sur $r\ge 1$ avec
\setcounter{equation}{0}
\begin{equation}
2^r \le x < 2^{r+1}\qquad\mathrm{et}\qquad y\ge 2.
\end{equation}
On choisit la chaîne
$$
\Cc(x,y) = \Big|\begin{array}{rl}
1-2 &\mathrm{si}\ r=1\\
4-1-2 &\mathrm{si}\ r=2
\end{array}
$$
et on vérifie qu'elle convient quand $x<8$. On suppose maintenant que $r\ge 3$, que toutes les chaînes $\Cc(x,y)$ avec $2\le x<2^r$ et $y\ge 2$ ont été construites et que (3.1) est vérifiée. On a en particulier $j\ge k\ge 1$.

Pour tout $x\ge 4$ on définit $\Cc^*(x,p_1)$ comme étant la chaîne constituée de tous les entiers $2^m$ vérifiant $2\le 2^m\le x$, par ordre décroissant.

On définit les deux suites finies d'entiers de $\Cs(x,y)$ suivantes :
$$
\Dd=p_1\Cc^*\Big(\frac{x}{p_1},p_1\Big)-p_2 
\Cc\Big(\frac{x}{p_2},p_2\Big)-p_3\Cc\Big(\frac{x}{p_3},p_3\Big)-\cdots-p_{k-1}\Cc\Big(\frac{x}{p_{k-1}},p_{k-1}\Big)
$$
et
$$
\Ee = \Big|
\begin{array}{l l}
p_k\Cc\Big(\dfrac{x}{p_k},p_k\Big) &\mathrm{si}\ k\ge 2\\
\emptyset &\mathrm{si}\ k=1.
\end{array}
$$

(Si $k\le 2$, on a donc $\Dd=p_1\Cc^*(x/p_1,p_1)$).

Si $k=1$, $\Dd$ est une chaîne.

Supposons à présent $k\ge 2$. On a alors $x\ge 18$ d'où
$$
j(x/p_2,p_2) = j(x/3,3) = \pi(3)=2.
$$
Le premier entier de $p_2\,\Cc(x/p_2,p_2)$ est donc égal à $12$, qui est un multiple de $4$, qui est le dernier entier de $p_1\Cc^*(x/p_1,p_1)$. De même pour $3\le \ell \le k$, en utilisant que $p_\ell \le p_k \le \sqrt{x/2}$, on~a
$$
j(x/p_\ell,p_\ell) = \pi (\min(p_\ell,x/2p_\ell))=\pi(p_\ell)=\ell. 
$$
Le premier entier de $p_\ell\,\Cc(x/p_\ell,p_\ell)$ est donc égal à $2p_{\ell-1}p_\ell$, qui est un multiple de $ép_{\ell-1}$, qui est le dernier entier de $p_{\ell-1}\Cc(x/p_{\ell-1},p_{\ell-1})$. Enfin par hypothèse de récurrence, tous les $p_\ell\,\Cc(x/p_\ell,p_\ell)$ sont des chaînes.

Tout cela prouve que
\begin{equation}
\Dd,\Ee\ \mathrm{et}\ \Dd -\Ee\ \textrm{sont\ des\ chaînes\ d'entiers\ de\ }\Cs(x,y).
\end{equation}

On définit maintenant l'inverse de $\Cc(x,y)$ par
$$
\inv\Cc(x,y) : \left|
\begin{array}{lll}
2-\Ee-1-\Dd &\mathrm{si} &j=k\\
2-1-\Dd-\Ee &\mathrm{si} &j=k+1\\
2-\Dd-\Ee -1-2p_{j-1} &\mathrm{si} &j\ge k+2.
\end{array}\right.
$$
On voit facilement que $\Cc(x,y)$ est une suite finie d'entiers de $\Cs(x,y)$ qui s'achève en 2. En  utilisant (3.2) et en discutant suivant les différents cas, on vérifie que $\Cc(x,y)$ est également une chaîne qui commence en $2p_{j-1}$. On a de plus
$$
\CA(x/2,y) = \{1\} \bigsqcup\ \bigsqcup_{p\le \min(y,\sqrt{x/2})} p\CA(x/2p,p).
$$
En réutilisant l'hypothèse de récurrence on voit donc  que $\Cc(x,y) \supset \CA(x/2,y)$. Cela achève la preuve de la proposition~3.1, et partant celle du théorème~1.1.

\section{La fonction $A(x,y,z,t)$}
 
 Pour $x$, $y$, $z$ et $t$ des réels positifs, on note $\CA(x,y,z,t):=\{n\le x : P(n)\le y,\ P^-(n)>z$ et $S(n)\le xt\}$ et $A(x,y,z,t) = |\CA(x,y,z,t)|$. On a (lemme~2 de \cite{Sai4}).
 
 \vspace{2mm}
 
 {\monlem Pour tout $(x,y,z,t)\in \R^{+*} \times [1,+\infty[^3$, on a}
 $$
 A(x,y,z,t)=\1_{[1,+\infty[}(x) + \sum_{z<p\le \min(y,\sqrt{xt})} A(x/p,p,z,t).
 $$

Par ailleurs on a le

 \vspace{2mm}
 
 \setcounter{equation}{0}
 
 {\monlem Soient $z_0\ge 1$ et $0<\varepsilon<1/2$. Sous les conditions
\begin{equation}
x\ge y \ge x^{1/2+\varepsilon}\ge 2\,,\ 1\le t\le x\ \mathrm{et}\ 1\le z \le z_0
\end{equation}
on a
\begin{equation}
A(x,y,z,t) \newsmile{\!}_{z_0,\varepsilon}\ x\ \frac{\log 2t}{\log x}\,.
\end{equation}}

\textit{Démonstration :} La majoration résulte du lemme 4 de \cite{Sai4}.  Pour la minoration, on a pour tout $t_0\ge 1$
$$
A(x,y,z,t) \ge A(x/t_0,y,z,tt_0)\,.
$$
Il suffit donc de prouver la minoration de (4.2) avec la condition supplémentaire $t\ge t_0$ où $t_0$ est une constante fixée. Or justement il résulte du lemme~1 de \cite{Sai5} qu'il existe une constante $t_0\ge 2$ telle que
\begin{equation}
\left[\begin{split}
A(x,+\infty,z,t) \gg_{z_0}\ x\ \frac{\log t}{\log x}\\
\mathrm{sous\ les\ conditions}\hbox to 4cm{}\\
x\ge t\ge t_0\ \mathrm{et}\ 1\le z\le z_0\,.
\end{split}
\right.
\end{equation}

On supposera donc dorénavant $t\ge t_0\ge 2$. Si $y>\sqrt{xt}$, on a en utilisant le lemme~4.1 et (4.3) que sous les conditions~(4.1)
$$
A(x,y,z,t) = A(x,+\infty,z,t) \gg_{z_0} \ x\ \frac{\log t}{\log x}\,.
$$

Si $y\le \sqrt{xt}$, il suffit de montrer la minoration de (4.2) pour $y=x^{1/2+\varepsilon}$. Or en utilisant le lemme~4.1 et (4.3) il existe une constante $x_0(\varepsilon)$ telle que sous (4.1) et la condition supplémentaire $x\ge x_0(\varepsilon)$ on~a
$$
\begin{array}{l}
A(x,x^{1/2+\varepsilon},z,t) \ge A(x,x^{1/2+\varepsilon},z,t) -A(x,\sqrt{x},z,t)\\ \\
= \displaystyle\sum_{\sqrt{x}<p\le x^{1/2+\varepsilon}} A(x/p,p,z,t) = \sum_{\sqrt{x}<p\le x^{1/2+\varepsilon}} A(x/p,+\infty,z,t)\\ \\
\gg_{z_0} \ x\ \frac{\log t}{\log x} \displaystyle\sum_{\sqrt{x}<p\le x^{1/2+\varepsilon}} \frac{1}{p} \newsmile{\!}_\varepsilon \ x\ \frac{\log t}{\log x}\,.
\end{array}
$$
On conclut la preuve de la minoration en remarquant qu'elle est banale si $x<x_0(\varepsilon)$.

\vspace{2mm}

{\monlem Il existe un entier positif $\gamma$ et un réel positif $x_0$ tels que pour tous réels $x$ et $q$ vérifiant 
\begin{equation}
x\ge x_0\quad \mathrm{et}\quad x^{1/5}<q\le x^{1/4},
\end{equation}
on a
$$
\Big|\Big\{\frac{\sqrt{x}}{q} <a \le \sqrt{\frac{x}{q}} : P(a) \le q\ \mathrm{et}\ S(a) \le \frac{2^\gamma\sqrt{x}}{q}\Big\}\Big| \ge \frac{\sqrt{x}}{2^\gamma q\log x}.
$$}

\noindent\textbf{Démonstration.} D'après le lemme $4.2$, il existe un réel $E\ge 2$ tel que
\begin{equation}
\left|
\begin{array}{l}
\mathrm{pour\ tous\ }X,\ q\ \mathrm{et}\ \alpha\ \textrm{vérifiant}\\
\vspace{3mm}
X\ge 2,\ X^{4/7}<q\le X,\ 1\le 2^\alpha \le X,\ \mathrm{on\ a}\\
\vspace{2mm}
\dfrac{1}{E} \le \dfrac{|\{a\le X:P(a)\le q\ \mathrm{et\ }S(a)\le 2^\alpha X\}|}{X(1+\alpha)/\log X} \le E
\end{array}
\right.
\end{equation}

On choisit un tel réel $E$. On choisit alors aussi un entier positif $\gamma$ tel que
\begin{equation}
2^\gamma \ge 4E^2 (1+\gamma).
\end{equation}
On choisit un réel $x_0$ tel que sous les conditions (4.4), on a
\begin{equation}
\sqrt{q} \ge 2^\gamma
\end{equation}
ce qui avec (4.6) entraîne
\begin{equation}
\sqrt{q} \ge 4E^2(1+8).
\end{equation}
Quitte à augmenter la valeur de $x_0$, on a alors en utilisant successivement (4.8) et (4.6), (4.5) et (4.4) et (4.8), et enfin de nouveau (4.5) et (4.4)
$$
\begin{array}{l}
\Big|\Big\{ a \le \sqrt{\dfrac{x}{q}}
 : P(a) \le q,\ S(a) \le \dfrac{2^\gamma\sqrt{x}}{q} \Big\}\Big|\\
 \ge \Big|\Big\{a\le 4E^2(1+\gamma) \dfrac{\sqrt{x}}{q} : P(a) \le q,\ S(a) \le 4E^2(1+\gamma) \dfrac{\sqrt{x}}{q} \Big\}\Big|\\
 \ge 8E(1+\gamma) \dfrac{\sqrt{x}}{q\log x} \\
 \ge 2 \Big|\Big\{ a \le \dfrac{\sqrt{x}}{q} : P(a) \le q,\ S(a) \le \dfrac{2^\gamma\sqrt{x}}{q} \Big\}\Big|
\end{array}
 $$
D'où avec (4.7), (4.5) et (4.6)
$$
\begin{array}{l}
\Big|\Big\{\dfrac{\sqrt{x}}{q} <a \le \sqrt{\dfrac{x}{q}} : P(a) \le q\ \mathrm{et\ } S(a) \le \dfrac{2^\gamma\sqrt{x}}{q}\Big\}\Big|\\
\ge \dfrac{1}{2} \Big|\Big\{ a \le \sqrt{\dfrac{x}{q}} : P(a) \le q \ \mathrm{et\ } S(a) \le \dfrac{2^\gamma\sqrt{x}}{q} \Big\}\Big|\\
\ge \dfrac{1}{2}\Big|\Big\{ a \le \dfrac{2^\gamma\sqrt{x}}{q} : P(a) \le q\ \mathrm{et\ } S(a) \le \dfrac{2^\gamma\sqrt{x}}{q}\Big\}\Big|\\
\ge \dfrac{1}{E} \dfrac{2^\gamma\sqrt{x}}{q\log x} \ge \dfrac{\sqrt{x}}{2^\gamma q\log x}.
\end{array}
$$

\section{Lemme de Gronwall Bidimensionnel Discret Arithmétique}

Rappelons que l'on note
$$
A(x,y,z) = A(x,y,z,1).
$$
\setcounter{equation}{0}

{\monlem (lemme GBDA)

Soient $h:]0,+\infty[\times [2,+\infty[\rg \R$ une fonction, $z\ge 1$ et $M\ge 1$ des réels tels que l'on a pour tous $x>0$ et $y\ge2$
\begin{equation}
 h(x,y) \ge 1\!\!1_{[1+\infty[}(x) + \sum_{z<p\le \min (y,\sqrt{x/M})} h(x/p,p).
\end{equation}
On a alors pour tous $x>0$ et $y\ge 2$
\begin{equation}
h(x,y) \ge A(x/M,y,z).
\end{equation}}

\noindent\textbf{Démonstration.} Par le lemme 4.1 avec $t=1$, on a l'identité de Buchstab
\begin{equation}
A(x,y,z) = 1\!\!1_{[1,+\infty[}(x) + \sum_{z<p\le \min(y,\sqrt{x})}A(x/p,p,z)
\end{equation}
pour tous $x>0$, $y\ge 2$ et $z\ge 1$.

Pour $0<x<M$, la somme en $p$ dans (5.1) est vide. On a donc $h(x,y)\ge 0=A(x/M,y,z)$ et (5.2) est vérifiée.

On procède à présent par récurrence sur $k\ge 0$ en prouvant que l'inégalité (5.2) est vérifiée pour tous les couples $(x,y)$ tels que $y\ge 2$~et
$$
0<x<M2^k \leqno(I_k)
$$
On l'a montrée quand $k=0$. On suppose à présent l'assertion vraie pour $k\ge 0$ et l'encadrement $(I_{k+1})$ vérifié.

En utilisant successivement (5.1), la minoration $M\ge 1$ et l'hypothèse de récurrence, et enfin (5.3), on obtient 
$$
\begin{array}{ll}
 h(x,y) &\ge 1\!\!1_{[1,+\infty[}(x) +\displaystyle\sum_{z<p\le \min(y,\sqrt{x/M})}h(x/p,p)\\
 &\ge 1\!\!1_{[1,+\infty[} (x/M) + \displaystyle\sum_{z<p\le\min(y,\sqrt{x/M})}A(x/Mp,p,z)\\
&= A(x/M,y,z),
\end{array}
$$
et (5.2) est vérifiée.

\section{Un lemme combinatoire}

Le lemme suivant est élémentaire. Il est déjà présent en  \cite{Sai2} comme le Lemme 5.

\vspace{2mm}

{\monlem Soit $\Ee$ un ensemble non vide partitionné en $\Ee= \coprod\limits_{j=1}^{k} \Ee_j$. On note $n=|\Ee|$.

Les conditions suivantes sont équivalentes.

$(i)$ On peut ordonner les éléments de $\Ee$ sous la forme $e_1,e_2,\ldots,e_n$ avec $\alpha_i \neq \alpha_{i+1}$ quand $e_i\in \Ee_{\alpha_i}$ pour tout $i$ tel que $1\le i<n$.}

\vspace{2mm}

$(ii)$ $\max\limits_{1\le j\le k}|\Ee_j| \le (n+1)/2$.

\vspace{2mm}

Nous donnons cette fois-ci la

\vspace{2mm}

\noindent\textbf{Preuve.} On suppose que
\setcounter{equation}{0}
\begin{equation}
|\Ee_1| \ge |\Ee_2| \ge \cdots \ge |\Ee_k|.
\end{equation}

$(i) \Rightarrow (ii)$

Notons $p$ et $q$ tels que 
$$
p=|\Ee_1| \qquad \mathrm{et} \qquad p+q=n.
$$
On suppose $(i)$. Alors l'application
$$
\begin{array}{ccl}
\Ee_1 \ba \{e_n\} &\lgr &\Ee\ba \Ee_1\\
e_i &\longmapsto &e_{i+1}
\end{array}
$$
est une injection. Donc
$$
p-1 \le |\Ee_1\ba \{e_n\}| \le |\Ee\ba \Ee_1|=q.
$$
D'où en rappelant que $q=n-p$, $p\le (n+1)/2$.

\vspace{2mm
}
$(ii) \Rightarrow (i)$

On procède par récurrence sur $n\ge 1$.

Les cas $n=1$ et $n=2$ sont faciles. Supposons $n\ge 3$ et l'implication vérifiée pour $n-2$. Comme $(n+1)/2 <n$, on a $k\ge 2$. Choisissons $f\in \Ee_1$ et $g\in \Ee_2$. On a $|\Ee_3| \le n/3 \le (n-1)/2$ car $n\ge3$. On a donc avec l'hypothèse $(ii)$ et (6.1),
$$
\max(|\Ee_1|-1, |\Ee_2| -1, \max_{3\le j\le k} |\Ee_j|) \le (n-1)/2.
$$
D'après l'hypothèse de récurrence pour $n-2$, on peut ordonner $\Ee\ba\{f,g\}$ en une suite $e_3,e_4,\ldots,e_n$ avec $\alpha_i \neq \alpha_{i+1}$ quand $e_i\in \Ee_{\alpha_i}$. En complétant avec $e_1$ et $e_2$ qui sont $f$ et $g$ tels que $e_2$ et $e_3$ ne sont pas dans le même $\Ee_j$, on obtient une suite $e_1,e_2,e_3,\ldots,e_n$ qui convient.

\section{Réserve de connecteurs}

Dans le cas de la preuve de la Proposition de \cite{Sai2}, on vérifie rapidement que l'ensemble $\Nn(x)$ (voir deuxième moitié de la page~241 de \cite{Sai2}) est suffisamment gros pour y puiser les connecteurs dont on a besoin pour la construction de la longue $x$--haute chaîne. L'objet du lemme suivant est d'effectuer la vérification analogue. La nécessité de choisir ici les connecteurs dans $\CA(x)$ augmente la longueur de cette vérification.

\vspace{2mm}

{\monlem Il existe un réel $x_0$ et un entier $\gamma\ge 0$ tels que pour tout $x\ge x_0$, pour tout ensemble $\CPp=\{p_{j_0}<p_{j_0+1}<\cdots<p_{k-1}\}$ non vide, qui d'une part est contenu dans l'intervalle $[13,\sqrt{x}/2^{\gamma+13}]$, et d'autre part contient tous les nombres premiers de l'intervalle $[13,p_{k-1}]$ à l'exception d'un ensemble $\Ee$ de nombres premiers de cardinal~$\le 2$, il existe une famille d'entiers $(q_ja_j)_{j_0\le j<k}$ deux à deux distincts et vérifiant pour tout $(p,q,a):=(p_j,q_j,a_j)$}
\setcounter{equation}{0}
\begin{equation}
\left[
\begin{array}{l}
q\ est\ premier,\ P(a)\le q,\ q\notin \Ee,\ 13\le q\le p,\\
\vspace{2mm}
\sqrt{x}/q<a\le \sqrt{x/q},\ 8pqS(a)\le x,
\end{array}
\right.
\end{equation}
\begin{equation}
2^{13} pq^2\le x
\end{equation}
\begin{equation}
q_j\not= q_{j+1}\ pour\ j_0\le j\le k-2.
\end{equation}

\noindent\textbf{Démonstration.} On supposera toujours implicitement $x$ suffisamment grand.

\vspace{2mm}

\textbf{Etape 1. Connecteurs intérieurs}

On suppose ici que $p=p_j$ vérifie
$$
13 \le p\le (x/2^{13})^{1/3}.
$$
On choisit alors $q=p$ et $a=2^{m(q)}$ tels que
$$
\frac{\sqrt{x}}{q} < a \le \sqrt{\frac{x}{q}}.
$$
On a alors $pq^2 \le p^3\le x/2^{13}$, d'où (7.2). On vérifie en fait que toutes les  conditions du lemme~7.1 sont satisfaites, la condition (7.3) uniquement pour les $j\le k-2$ tels que $p_j\le (x/2^{13})^{1/3}$.

\vspace{2mm}

\textbf{Etape 2. Connecteurs extérieurs}

On se donne un entier positif $\gamma$ quelconque pour l'instant et on suppose ici que pour $p=p_j$
$$
(x/2^{13})^{1/3} < p \le \sqrt{x}/2^{\gamma+13}
$$
On recouvre l'intervalle $](x/2^{13})^{1/3},\ \sqrt{x}/2^{\gamma+13}]$ par les intervalles
$$
I_r=\Big]\frac{\sqrt{x}}{2^{\gamma+13+r+1}},\ 
\frac{\sqrt{x}}{2^{\gamma+13+r}}\Big]
$$
avec $0\le r\le R$ où $R=R(x)$ est l'entier maximum tel que
$$
(x/2^{13})^{1/3} < \sqrt{x} /2^{\gamma+13+R}
$$
On note
$$
q_r(x): = x^{\frac{1}{4}(1-\frac{1}{9}\sum\limits_{j=1}^r\frac{1}{j^2})}
$$
En utilisant le théorème des nombres premiers, on a pour tout $r\in[O,R]$
$$
\begin{array}{l}
2 \displaystyle \sum_{q_{r+1}<q\le q_r\atop q\notin \Ee} \dfrac{1}{q} \ge \log \Big(\dfrac{\log q_r}{\log q_{r+1}}\Big)\\
=\log\Big(1- \dfrac{1}{9} \displaystyle\sum_{j=1}^r \dfrac{1}{j^2}\Big) -\log \Big(1 - \dfrac{1}{9} \displaystyle\sum_{j=1}^{r+1} \dfrac{1}{j^2}\Big)\\
\ge \dfrac{1}{9(r+1)^2}.
\end{array}
$$
On en déduit avec le lemme 4.3 qu'avec maintenant un entier $\gamma\ge 0$ ad hoc, on a
$$
\begin{array}{l}
\displaystyle\sum_{q_{r+1}<q\le q_r\atop q\notin \Ee} \Big|\Big\{\dfrac{\sqrt{x}}{q} < a \le \sqrt{\dfrac{x}{q}} : P(a) \le q,\ S(a) \le \dfrac{2^\gamma\sqrt{x}}{q},\ a\not= 2^{m(q)}\Big\}\Big|\\
\ge \dfrac{\sqrt{x}}{2^{\gamma+1}\log x} \displaystyle\sum_{q_{r+1}<q\le q_r\atop q\notin \Ee} \dfrac{1}{q} \ge \dfrac{1}{9\cdot 2^{\gamma+2}(r+1)^2} \dfrac{\sqrt{x}}{\log x}\\
\ge \dfrac{1}{2^{\gamma+12+r}} \dfrac{\sqrt{x}}{\log x}\\
\ge \pi \Big(\dfrac{\sqrt{x}}{2^{\gamma+13+r}}\Big) - \pi \Big(\dfrac{\sqrt{x}}{2^{\gamma+13+r+1}}\Big).
\end{array}
$$
par le théorème des nombres premiers.

Remarquons également que
$$
q_{r+1} <q \le q_r \Rightarrow q\le x^{1/4} < (x/2^{13})^{1/3} <p \le \sqrt{x}/2^{\gamma+13}.
$$
Donc pour tout $r\in [O,R]$, il existe une application injective
$$
\begin{array}{rl}
\varphi_r : \CPp \cap I_r &\longrightarrow \{n\ge 1 : q_{r+1} <P(n) \le q_r\}\\
p &\longmapsto qa = q(p)a(p)
\end{array}
$$
telle que les conditions (7.1) sont vérifiées.

On a
$$
p\le \sqrt{x}/2^{\gamma+13}\le \sqrt{x}/2^{13}
$$
et $q\le x^{1/4}$. Cela montre que la condition (7.2) est également vérifiée.

Soit $0\le r\le R$. On a $2^r \le 2^R \le x^{1/6}$ et $q_{r+1}>x^{1/5}$. On en déduit que pour tout $q\in ]q_{r+1},q_r]$,
$$
\begin{array}{l}
\Big|\Big\{\dfrac{\sqrt{x}}{q} <a\le \sqrt{\dfrac{x}{q}} : P(a) \le q,\ S(a) \le \dfrac{2^\gamma\sqrt{x}}{q},\ a\not=2^{m(q)} \Big\}\Big|\\
\le \dfrac{2^\gamma\sqrt{x}}{q}\\
\le \dfrac{1}{2}\Big(\pi \Big(\dfrac{\sqrt{x}}{2^{\gamma+13+r}}\Big) - \pi \Big(\dfrac{\sqrt{x}}{2^{\gamma+13+r+1}}\Big)+1\Big).
\end{array}
$$
de nouveau par le théorème des nombres premiers.

D'après le lemme 6.1, on peut donc choisir l'application injective $\varphi_r$ de telle sorte que l'on a $q(p)\not=q(p^-)$. Ainsi la condition (7.3) est satisfaite pour les $q$ appartenant à un même intervalle $]q_{r+1},q_r]$.

Comme les intervalles $]q_{r+1},q_r]$ sont disjoints cela entraîne que l'on a aussi
$$
q(p) \not=q(p^-)
$$
pour tous les $p=p_j$ vérifiant $(x/2^{13})^{1/3}<p\le \sqrt{x}/2^{\gamma+13}$. Pour le plus petit $p>(x/2^{13})^{1/3}$, on a $q(p) \le x^{1/4}<q(p^-)$. Tout cela montre que la condition (7.3) est satisfaite.

Enfin la condition $a\not=2^{m(q)}$ assure que les $qa$ choisis pour $p>(x/2^{13})^{1/3}$ sont différents de ceux pour $p\le (x/2^{13})^{1/3}$. Cela achève de prouver que les entiers $(q_ja_j)_{j_0\le j<k}$ que l'on a choisi sont deux à deux distincts. Cela termine la preuve du lemme~7.1.

\section{Construction de chaînes d'entiers de $\CA^*(x,y)$}

L'objet du chapitre 8 de \cite{Sai2} a été d'établir une minoration de $f^*(x)$ par une fonction de comptage d'entiers à diviseurs denses, minoration qui s'est révélée par la suite donner le bon ordre de grandeur de $f^*(x)$ (on a $f^*(x)\asymp f(x) \asymp x/\log x$, cf. le théorème~3 de \cite{Sai4}). On a construit pour cela des chaînes  d'entiers de $\{\sqrt{x} <n \le x : P(n) \le y\}$ qui ont la structure (1.3), ce qui a permis d'établir une inéquation de Buchstab pour une fonction $h^*(x,y)$ qui minore $f^*(x,y)$, que l'on a itérée.

Pour montrer que
$$
f_a^*(x) \gg x/\log x
$$
on suit la même procédure. C'est la fonction $h_a^*(x,y)$ qui sera introduite à la formule (8.5), qui va jouer le rôle de $h^*(x,y)$ dans \cite{Sai2}. On va en fait suivre de très près la démarche de \cite{Sai2}. On espère cependant que la rédaction est meilleure ici.

On va construire des chaînes d'entiers de
$$
\CA^*(x,y) := \{n>\sqrt{x} : P(n)\le \min(y,\sqrt{x}/27),\ S(n)\le x\}
$$
qui vont relier un entier $A$ à un entier $B$, en collant trois morceaux
$$
\Cc_{A,B}^* : \Cc_{A,A'}^* - \Cc_{A',B'}^* - \Cc_{B',B}^*
$$
comme dans \cite{Sai2} (voir le milieu de la page~238). C'est la chaîne centrale $\Cc_{A',B'}^*$ qui a une structure de la forme (1.3) et qui va permettre d'établir l'inégalité de Buchstab de $h_a^*(x,y)$. Les entiers $n$ de $\Cc_{A',B'}^*$ doivent pour la plupart vérifier $n>\sqrt{xP(n)}$, ce qui n'est pas nécessairement le cas des entiers $A$ et $B$. C'est la raison pour laquelle, comme dans \cite{Sai2} $(c2)$ construction de $\Cc_{A,A'}$ et $\Cc_{B',B}$ de la page~239), on est amené d'abord à construire les chaînes $\Cc_{A,A'}$ et $\Cc_{B',B}$ qui permettent de relier $A$ à $A'$ et $B$ à $B'$ avec $A'$ et $B'$ vérifiant $A'>\sqrt{xP(A')}$ et $B'>\sqrt{xP(B')}$.

Le cas $y=3$ est particulier et la construction de la chaîne $\Cc_{A,B}^*$ dans ce cas est faite au lemme~8.5.

Pour le cas générique $y\ge5$, on construit les chaînes $\Cc_{A,A'}^*$ et $\Cc_{B',B}^*$ au lemme~8.6. L'essentiel du travail réside dans la construction de $\Cc_{A',B'}^*$ du lemme~8.7. Les lemmes de 8.1 à 8.4 constituent un travail préparatoire.

\vspace{2mm}

{\monlem Soient $q\in \{2,3\}$, $x$ et $y$ deux réels tels que $y\ge 3$ et
\setcounter{equation}{0}
\begin{equation}
x\ge 3^{10}.
\end{equation}

Soit $A$ un entier de $\CA^*(x,y)$ tel que
\begin{equation}
q<P(A) \le\sqrt{x}/27
\end{equation}
et
$$
P^-(A)\ge q.
$$

Il existe alors une chaîne $\Cc(x,y) : A=a_1-a_2 -\cdots -a_s$ d'au moins deux entiers de $\CA^*(x,y)$ telle que}
$$
\begin{array}{l}
P(a_j)=P(A)\ pour\ tout\ j\ \textit{vérifiant}\ 1\le j <s,\\
P(a_s)=q,\\
P^-(a_j) \ge q\ pour\ tout\ j\ \textit{vérifiant}\ 1\le j\le s,
\end{array}
$$
\begin{equation}
\left[
\begin{array}{ll}
&a_{s-1}>\sqrt{xP(a_{s-1})},\\
&P(a_{s-1}/P(a_{s-1})) \le \sqrt{x/P(a_{s-1})}/27,\\
et\\
&a_s >3\sqrt{x}
\end{array}
\right.
\end{equation}

\noindent\textbf{Remarque.} Quand $q=2$, les minorations de l'énoncé du type $P^-(n)\ge q$ sont vides.

\vspace{2mm}

\noindent\textbf{Démonstration.} Les hypothèses sur l'entier $A$ permettent de l'écrire sous la forme
$$
A=p_1p_2\cdots p_kq^\alpha
$$
avec $k\ge 1$ et
$$
P(A) =p_1\ge p_2\ge \cdots \ge p_k>q.
$$
On note $A'=p_1p_2\cdots p_k q^{\alpha+\nu_1}$ l'entier défini par
\begin{equation}
x/q^2 <A' \le x/q.
\end{equation}
Comme $S(A) \le x$, on a $A\le x/q$. On en déduit que $\nu_1\ge 0$. Cette majoration de $A'$ combinée avec l'hypothèse $A\in \CA^*(x,y)$ assurent que $A'\in \CA^*(x,y)$. On définit
$$
A'' := A'/p_k.
$$
La minoration de (8.4), la majoration $p_k \le \sqrt{x}/27$ et  le fait que $A'\in A^*(x,y)$ assurent que $A'' \in \CA^*(x,y)$.

Si $k=1$, on choisit la chaîne
$$
\Cc(x,y)= A-A'-A''.
$$

\noindent\textbf{Remarque.} Ce dernier schéma est à comprendre de la manière suivante. Si $A'\not= A$, la chaîne $\Cc(x,y)$ est constituée des trois entiers $A$, $A'$ et $A''$. Si au contraire $A'=A$, la chaîne $\Cc(x,y)$ est constituée uniquement des deux entiers $A$ et $A''$. Autrement dit la liaison $A-A'$ est suivant les cas une concaténation ou un collage. On utilisera  par la suite plusieurs fois cette notation --- qui a cette signification différente suivant les cas, sans plus le faire remarquer.

Si $k\ge 2$, on prolonge la chaîne $A-A'-A''$ vers la droite en répétant cette paire d'opérations, avec
$$
A'' = p_1p_2 \cdots p_{k-1}q^{\alpha+\nu_1}
$$
dans le rôle de $A$. La nouvelle chaîne est donc $A-A'-A''-A'''-A''''$ avec $A'''=p_1p_2\cdots p_{k-1}q^{\alpha+\nu_1+\nu_2}$ et $A''''=p_1p_2 \cdots p_{k-2}q^{\alpha+\nu_1+\nu_2}$. Remarquons que cette fois--ci on a $\nu_2>0$ et donc que $A''-A'''$ est toujours une concaténation. En itérant, on effectue au final $k$ fois cette paire d'opérations qui permet à chaque fois de diminuer de un le nombre, compté avec multiplicité, de facteurs premiers différent de $q$. On obtient ainsi une chaîne $\Cc(x,y)$ d'entiers de $\CA^*(x,y)$, ayant les propriétés demandées pour les plus grand et plus petit facteurs premiers des entiers dont elle est constituée. Notons
$$
\Cc(x,y) = A =a_1-a_2 -\cdots - a_{s-1}-a_s
$$
cette chaîne.

Montrons à présent les inégalités (8.3). On a $P(a_{s-1})=P(A)$.

Comme dans la minoration de (8.4), on a $a_{s-1}> x/q^2$. Avec (8.1) et la majoration de (8.2), on a donc
$$
\begin{array}{c}
a_{s-1} >x/9,\\
x\ge 3^{10}
\end{array}
$$
et
$$
P(a_{s-1}) \le \sqrt{x} /27.
$$
Cela permet de montrer que
$$\begin{array}{c}
a_{s-1} >\sqrt{xP(a_{s-1})},\ a_s=a_{s-1}/P(a_{s-1})> 3\sqrt{x}\\
\vspace{2mm}
\mathrm{et}\ P(a_{s-1}/P(a_{s-1})) = q \le 3 \le \sqrt{x/P(a_{s-1})}/27.
\end{array}
$$

\vspace{2mm}

{\monlem Pour tous réels $x\ge 3^{10}$ et $y$ et tout entier $A$ de $\CA^*(x,y)$ tel que $P(A)> 3$, il existe une chaîne $A=a_1-a_2-\cdots -a_s$ d'éléments de $\CA^*(x,y)$ telle que $P(a_j)=P(A)$ pour tout $j$ vérifiant $1\le j<s$, et $P(a_s)=3$.}

\vspace{2mm}

\noindent\textbf{Démonstration.} Si  $P^-(A)\ge 3$, on pose $r=1$. Si $P^-(A)=2$, on commence par construire une chaîne $A=a_1-a_2-\cdots -a_r$ d'entiers de $\CA^*(x,y)$ avec $r\ge 2$ telle que

\vspace{2mm}

\begin{description}
\item[*] Pour tout $j$ vérifiant $2\le j\le r$, $a_j/a_{j-1}$ est une puissance de 3 ou  $a_{j-1}/a_j$ est une puissance de 2

\item[*] Pour tout $j$ vérifiant $1\le j\le r-1$, $P^-(a_j)=2$

\item[*] $P^-(a_r)=3$.
\end{description}

\vspace{2mm}

Comme $P(a_r)=P(A)$ et $A\in \CA^*(x,y)$, on a
$$
3 < P(a_r) \le \sqrt{x}/27.
$$
On a également $x\ge 3^{10}$.

On peut donc appliquer le lemme 8.1 avec $q=3$, ce qui nous fournit une chaîne $\Cc(x,y)$ d'entiers de $\CA^*(x,y)$ qui commence à $a_r$ 
$$
\Cc(x,y) : a_r - a_{r+1} -\cdots - a_s,\ \mathrm{avec}
$$
$$
\begin{array}{l}
P(a_j) =P(a_r)\ \mathrm{pour\ } r\le j <s\\
P(a_s) =3\\
P^-(a_j) \ge 3\ \mathrm{pour}\ r\le j\le s.
\end{array}
$$
Cela permet d'achever la preuve que la suite finie
$$
A=a_1-a_2 -\cdots - a_r - a_{r+1} -\cdots - a_s
$$
est une chaîne qui répond à la question.

\vspace{2mm}

{\monlem Pour tous réels $x\ge 3^{10}$ et $y\ge 3$, pour tous $A$ et $B$ éléments de $\CA^*(x,y)$ tels que $P(A)<P(B)$, il existe une chaîne formée d'entiers de $\CA^*(x,y)$, qui relie $A$ à $B$.}

\vspace{2mm}

\noindent\textbf{Démonstration.} Si $P(A)=2$, on choisit $\Cc_1(x,y)$ réduit au singleton $A$. Sinon on choisit la chaîne $\Cc_1(x,y)$ comme la chaîne $\Cc(x,y)$ du lemme 8.1 pour $q=2$, qui commence par $A$. Par ailleurs on choisit $\Cc_2(x,y)$ également comme la chaîne $\Cc(x,y)$ du lemme~8.1 pour $q=2$, mais qui commence cette fois par $B$. On obtient finalement une chaîne qui convient en concaténant ou en collant $\Cc_1(x,y)$ à l'inverse de $\Cc_2(x,y)$

\vspace{2mm}

{\monlem Soit $x\ge 3^8$.

Il existe alore deux entiers $A$ et $B$ de $\CA^*(x,3)$ tels que $P(A)=2$ et $P(B)=3$.}

\vspace{2mm}

\vspace{2mm}

\noindent\textbf{Démonstration.} On choisit $A$ comme la plus grande puissance de 2 qui est $\le x/6$. On choisit alors $B=3A$. Comme $x\ge 12$, $A\not=1$ et $S(B)=2B$. Plus précisément, $S(A) =2A<6A=S(B)\le x$. Comme $x\ge 3^8$, on a aussi
$$
2=P(A) < P(B)=3 \le \sqrt{x}/27.
$$

Le lemme 8.3 nous permet de définir $f_{a,A,B}^*(x,y)$ comme la longueur maximum d'une chaîne d'entiers de $\CA^*(x,y)$ reliant $A$ à $B$, dès que\break $(A,B)\in \CA^*(x,y)^2$, $P(A)<P(B)$, $x\ge 3^{10}$ et $y\ge3$. Avec le lemme~8.4, cela permet de définir également la fonction $h_a^*$ de $\R^{+*}\times [2,+\infty[$ dans $\N$~par
\begin{equation}
h_a^*(x,y)= \left|
\begin{array}{l}
{1\!\!1_{[1,+\infty[}(x)}\ \mathrm{si}\ [0<x<3^{10}\ \mathrm{ou}\ 2\le y<3]\\
\min\limits_{(A,B)\in \CA^{*2}(x,y),P(A)<P(B)} f_{a,A,B}^* (x,y),\ \mathrm{sinon}.
\end{array}
\right.
\end{equation}
Dans la suite  de ce chapitre, la notation $\Cc_{C,D}^*$ désignera toujours, parfois implicitement, une chaîne commençant en $C$ et s'achevant en~$D$.

\vspace{2mm}

{\monlem Soit $x\ge  3^{11}$.

Soient $A$ et $B$ deux entiers de $\CA^*(x,3)$ avec $P(A)=2$ et $P(B)=3$. Il existe alors une chaîne $\Cc_{A,B}^*$ d'entiers de $\CA^*(x,3)$ qui relie $A$ à $B$ avec}
$$
longueur\ \Cc_{A,B}^* \ge 1 + h_a^*(x/3,3).
$$

\noindent \textbf{Remarque.} A plusieurs occasions dans cette preuve, on utilisera implicitement une hypothèse de minoration de $x$, plus faible que $x\ge 3^{11}$.

\vspace{2mm}

\noindent\textbf{Démonstration.} On choisit
$$
B'= \left|
\begin{array}{c l}
B &\mathrm{si\ }B>\sqrt{3x}\\
2B &\mathrm{sinon.}
\end{array}
\right.
$$
Si $B'=2B$, on a $S(B')=4B \le 4\sqrt{3x}\le x$. On a de plus $P(B')=3\le\sqrt{x}/27$. Donc $B'\in \CA^*(x,3)$. On note
$$
\Cc_{B',B} : B'-B.
$$

Par ailleurs on note $r$ l'entier défini par $\sqrt{x} < 2^r \le 2\sqrt{x}$. On choisit alors un entier $A'$ et une chaîne de la forme
$$
\Cc_{A,A'}^* : A-2^r -2^r 3^{1+\varepsilon}:= A'
$$
avec $\varepsilon\in \{0,1\}$ choisi de telle sorte que $P(A'/3)\not= P(B'/3)$. On a $S(A')=2A'\le 2^2\cdot 3^2 \sqrt{x}\le x$. On a toujours $3\le \sqrt{x}/27$. Donc $\Cc_{A,A'}^*$ est une chaîne d'entiers de~$\CA^*(x,y)$.

On a $A'/3>\sqrt{x/3}$, $B'/3>\sqrt{x/3}$, $3\le\sqrt{x/3}/27$, $S(A'/3)\le x/3$, $S(B'/3)\le x/3$ et $P(A'/3)\not= P(B'/3)$. On a de plus $x/3\ge 3^{10}$ d'après l'hypothèse $x\ge 3^{11}$.

D'après le lemme 8.3 qui est utilisé pour définir $h_a^*$ à la formule (8.5),\break il existe une  chaîne $\Cc_{A'/3,B'/3}^*$ d'entiers de $\CA^*(x/3,3)$ et de longueur\break $f_{a,A',B'}^*(x/3,3)\ge h_a^*(x/3,3)$. De plus l'entier $A$ n'appartient pas à la chaîne $3\Cc_{A'/3,B'/3}^*$. Cela permet de conclure que
$$
\Cc_{A,B}^* =\Cc_{A,A'}^* -3\Cc_{A'/3,B'/3}^* - \Cc_{B',B}^*
$$
est une chaîne d'entiers de $\CA^*(x,3)$ qui relie $A$ à $B$ et dont la longueur est supérieure ou égale à $1+h_a^*(x/3,3)$.

Passons maintenant au cas générique $y\ge 5$. Comme on l'a déjà mentionné, on construit une chaîne $\Cc_{A,\CPp}^*$ de la forme
\begin{equation}
\Cc_{A,B}^* : \Cc_{A,A'}^* - \Cc_{A',B'}^* -\Cc_{B',B}^*.
\end{equation}

 \vspace{2mm}
 
 {\monlem \textrm{(Première partie : construction de $\Cc_{A,A'}^*$ et $\Cc_{B,B'}^*$)}
 
 Soit $\gamma$ un entier positif ou nul.
 
 Il existe alors un réel positif $x_1$ tel que pour tous réels $x$ et $y$ vérifiant
 $$
 x\ge x_1\ et\ y\ge 5,
 $$
 pour tous entiers $A$ et $B$ éléments de $\CA^*(x,y)$ vérifiant
 $$
 P(A) < P(B),
 $$
 il existe deux entiers $A'$ et $B'$ et deux chaînes disjointes  $\Cc_{A,A'}^*$ et $\Cc_{B,B'}^*$ d'entiers de $\CA^*(x,y)$ qui vérifient} 
 \begin{equation}
 il\ y\ a\ au\ plus\ une\ puissance\ de\ 2\ dans \ \Cc_{A,A'}^*\cap]\sqrt{x},2^8 \sqrt{x}]
 \end{equation}
 \begin{equation}
 il\ y\ a\ au\ plus\ une\ puissance\ de\ 2\ dans \ \Cc_{B,B'}^*\cap]\sqrt{x},2^8 \sqrt{x}]
 \end{equation}
 \begin{equation}
 P(A') \not= P(B')
 \end{equation}
 \begin{equation}
 n\in \Cc_{A,A'} \sqcup \Cc_{B,B'} \Rightarrow \left|
 \begin{array}{l}
 P(n)\in \{2,P(A')P(B')\}\\
 ou\ P(n)>\sqrt{x}/2^{\gamma+13}
 \end{array}
 \right.
 \end{equation}
 \begin{equation}
 \left|
 \begin{array}{c}
 les \ entiers\ A'\ B'\ sont \ les\ seuls\ entiers\ n\\
 de\ \Cc_{A,A'}^*\sqcup \Cc_{B,B'}^*\ qui\ \textit{vérifient}\\
 \textit{simultanément\ les \ trois\ propriétés\ suivantes } :\\
 3\le P(n) \le \sqrt{x}/2^{\gamma+13}\\
 n>\sqrt{xP(n)}\\
 P(n/P(n)) \le \sqrt{x/P(n)}/27.
 \end{array}
 \right.
 \end{equation}
 
 \vspace{2mm}
 
 \textbf{Remarques.} Dans la preuve qui suit
 
 \begin{description}
 \item[1)] On explicite la construction de $\Cc_{A,A'}^*$ et $\Cc_{B,B'}^*$ en distinguant différents cas.
 
 \item[2)] Dans chaque cas, il faut vérifier que les différentes propriétés de $\Cc_{A,A'}^*$ et $\Cc_{B,B'}^*$ demandées par ce lemme sont satisfaites. Ces vérifications ne posent pas de difficulté ; on les omettra.
 
 \item[3)] Quel que soit le cas étudié, les nombres premiers $p_o$ et $p_\infty$ seront définis par
 $$
 p_o := P(A') \ \mathrm{et}\ p_\infty := P(B').
 $$

 \item[4)] Au long de la preuve, le réel $x$ sera toujours supposé implicitement suffisamment grand.
 \end{description}
 
 \vspace{2mm}
 
\noindent\textbf{Démonstration. Etape 1. Construction de $\Cc_{A,A'}$}

\vspace{2mm}

\textbf{1\up{er} cas} $P(A)=2$

On choisit $p_o\in\{3,5\}$ tel que $p_o\not= P(B)$. On définit l'entier $\alpha$ par
$$
x/2^5 < 2^\alpha \le x/2^4
$$
on choisit alors la chaîne
$$
\Cc_{A,A'} : A-2^\alpha - 2^\alpha p_o := A'
$$

\vspace{2mm}

\textbf{2\up{ième} cas}
$$
3\le P(A) \le \sqrt{x}/2^{\gamma+13},\ A>\sqrt{xP(A)}\ \mathrm{et\ }P(A/P(A)) \le \sqrt{x/P(A)}/27.
$$

On choisit $A'=A$ et donc $\Cc_{A,A'}$ réduite au seul entier~$A$.

\vspace{2mm}

\textbf{3\up{ième} cas}
$$
3\le P(A) \le \sqrt{x}/2^{\gamma+13},\ \mathrm{et\ } [A\le\sqrt{xP(A)}\ \mathrm{ou\ }P(A/P(A)) > \sqrt{x/P(A)}/27].
$$

On note $A=a_1 -a_2 -\cdots - a_s=2^{\alpha'}$ la chaîne construite au lemme 8.1 avec $q=2$. On choisit $\ell$ l'entier minimum $\ge 2$ tel que
$$
a_\ell >\sqrt{xP(a_\ell)} \ \mathrm{et\ } P(a_\ell/P(a_\ell)) \le \sqrt{x/P(a_\ell)}/27.
$$
On sait qu'il en existe car l'entier  $s-1$ vérifie ces conditions. On choisit alors
$$
\Cc_{A,A'} : A= a_1-a_2 -\cdots- a_\ell=A'.
$$

\vspace{2mm}

\textbf{4\up{ième} cas} $P(A)>\sqrt{x}/2^{\gamma+13}$

On note 
$$
\Cc_{A,2^r} : A =a_1 - a_2 -\cdots -a_s := 2^{\alpha'}
$$
la chaîne construite au lemme 8.1 avec $q=2$. On définit comme au 1\up{er} cas l'entier $\alpha$ par
$$
x/2^5 < 2^\alpha \le x/2^4
$$ 
On choisit $p_o\in \{3,5\}$ tel que $p_o\not=P(B)$. On choisit enfin la chaîne
$$
\Cc_{A,A'} : \Cc_{A,2^{\alpha'}} -2^\alpha -2^\alpha p_o:= A'.
$$

\vspace{2mm}

\noindent\textbf{Etape 2. Construction de $\Cc_{B,B'}^*$}

On construit la chaîne $\Cc_{B,B'}^*$ de manière analogue. On différentie les mêmes quatre cas en remplaçant $A$ par $B$. Remarquons cependant que comme $P(B)\ge 3$, on n'est jamais dans le 1\up{er} cas pour $B$.

Dans les 2\up{ième} et 3\up{ième} cas pour $B$, on procède comme dans les 2\up{ième} et 3\up{ième} cas pour $A$, avec l'utilisation du lemme~8.1 avec $q=2$ dans le 3\up{ième} cas.

Dans le 4\up{ième} cas, on commence la construction de $\Cc_{B,B'}^*$ aussi de manière analogue. On note $\Cc_{B,2^{\beta'}}^*$ la chaîne fournie par le lemme~8.1 avec $q=2$.
D'après la dernière inégalité de (8.3), on a $2^{\beta'}>2\sqrt{x}$. Pour $\beta''$ l'entier défini par $2\sqrt{x}<2^{\beta''}\le 4\sqrt{x}$, on a donc $\beta''\le\beta'$. D'après (8.7), on peut choisir $\beta'''\in \{\beta'',\beta''-1\}$ tel que $2^{\beta'''}\notin \Cc_{A,A'}^*$. On choisit $p_\infty\in \{3,5\}$ tel que $p_\infty\not=p_o$. Notons
$\Cc_{B,2^{\beta'}}^*{\kern-0.5cm}'$\kern0.5cm 
la chaîne obtenue à partir de $\Cc_{B,2^{\beta'}}^*$ en supprimant son dernier entier $2^{\beta'}$. On choisit alors
$\Cc_{B,B'}^* : \Cc_{B,2^{\beta'}}^*{\kern-5mm}'\kern5mm - 2^{\beta'''}- 2^{\beta'''}p_\infty$. Pour vérifier que $\Cc_{B,B'}^* \subset \CA^*(x,y)$, on utilise en particulier que $2^{\beta'''}\ge 2^{\beta''-1}>\sqrt{x}$.

\vspace{2mm}

{\monlem (deuxième partie : construction de $\Cc_{A',B'}^*$)

Il existe un entier $\gamma\ge 0$ et un réel positif $x_2$ tels que pour tous réels $x$ et $y$ vérifiant
\begin{equation}
5 \le y\le \sqrt{x}/2^{\gamma+13}\ \mathrm{et\ }x\ge x_{2},
\end{equation}
pour tous entiers $A$ et $B$ éléments de $\CA^*(x,y)$ vérifiant $P(A)<P(B)$, il existe une chaîne $\Cc_{A,B}^*$ d'entiers de $\CA^*(x,y)$ qui relie $A$ à $B$ avec}
$$
longueur\ \Cc_{A,B}^* \ge 1 + \sum_{2<p\le y} h_a^*(x/p,p).
$$

\vspace{2mm}

\noindent\textbf{Remarque.} Dans toute cette preuve, on supposera toujours implicitement que $x$ est suffisamment grand.

\vspace{2mm}
 
\noindent\textbf{Démonstration. Etape 1. Structure de la chaîne $\Cc_{A,B}^*$}

\vspace{2mm}

On choisit un entier positif $\gamma$ qui convient pour le lemme 7.1. On choisit deux chaînes d'entiers $\Cc_{A,A'}^*$ et $\Cc_{B,B'}^*$ vérifiant les propriétés du lemme~8.6. On réutilise les notations
$$
p_o := P(A')\ \mathrm{et\ }p_\infty := P(B').
$$
On récrit aussi la propriété (8.9) par
\begin{equation}
p_o \not= p_\infty.
\end{equation}

On pose $\Cc_{B',B}^*=\inv \Cc_{B,B'}^*$. L'objet de ce lemme 8.7 est de compléter la construction du collage~(8.6) en construisant la partie centrale $\Cc_{A',B'}^*$. On va la construire sous la forme (1.3) avec parfois l'insertion de connecteurs extérieurs. On désignera par $p$ l'un des nombres premiers $p$ qui intervient dans (1.3). Autrement dit, on utilisera la lettre $p$ pour désigner un nombre premier générique vérifiant
$$
3\le p \le y.
$$
Les sous--chaînes $\Cc(x/p_j,p_j)$ de (1.3) seront notées ici $\Cc^*(x/p,p)$.

Notons $k:=\pi(y)-3$. Comme $y\ge 5$, on a
\begin{equation}
k\ge 0.
\end{equation}

Ordonnons les $p$. On a déjà défini $p_o$ et $p_\infty$. Si $k\ge 1$, on note également
$$
p_1 < p_2 <\cdots <p_k
$$
les nombres premiers $p$ restant de telle manière que
\begin{equation}
\left|
\begin{array}{cll}
\mathrm{l'application\ }\eta:\{0,1,2,\ldots,k\} \sqcup \{+\infty\} &\lgr &\{p:3 \le p\le y\}\\
j &\longmapsto &p_j
\end{array}
\right.
\end{equation}
est une bijection. C'est possible d'après (8.13).

On précise la structure (1.3) de $\Cc_{A',B'}^*$ par le schéma suivant :
%\eject

\begin{equation}
\unitlength=1cm
\begin{picture}(5,2)
\put(-2.6,1.6){$p_o\Cc^*(\frac{x}{p_o},p_o)$}
%1
\put(-1,1.2){\line(1,-2){0.5}}
\put(1,1.6){$p_1\Cc^*(\frac{x}{p_1},p_1)$}
%2
\put(1.2,0.2){\line(1,2){0.5}}
%3
\put(2,1.2){\line(1,-2){0.5}}
%4
\put(3.6,0.2){\line(1,2){0.3}}
\put(4.2,0.8){$\cdots$}\\
%5
\put(5,0.8){\line(1,-2){0.3}}
\put(6,1.6){$p\Cc^*(\frac{x}{p},p)$}
%6
\put(6,0.2){\line(1,2){0.5}}
%deuxième ligne
\put(-1,-3.4){\line(1,2){0.5}}
\put(-1.6,-2){$p^+\Cc^*(\frac{x}{p^+},p^+)$}
%8
\put(7,1.2){\line(1,-2){0.5}}
%10
\put(1.4,-3){\line(1,2){0.3}}
\put(0.6,-3.3){$\cdots$}\\
%9
\put(0,-2.4){\line(1,-2){0.3}}
\put(1.6,-2){$p_k\Cc^*(\frac{x}{p_k},p_k)$}
\put(5,-2){$p_\infty\Cc^*(\frac{x}{p_\infty},p_\infty)$}
%11
\put(3.2,-2.4){\line(1,-2){0.5}}
%12
\put(5,-3.4){\line(1,2){0.5}}
\put(0,-0.5){$q_oa_o$}
\put(2.6,-0.5){$q_1a_1$}
\put(5,-0.5){$q(p^-)a(p^-)$}
\put(-2,-4){$qa=q(p)a(p)$}
\put(4,-4){$q_k a_k$}
\put(-3,-4.2){\line(0,1){6.5}}
%\put(2,-2.2){$\Cc_{A',\CPp'}^*$}
\end{picture}
\end{equation}
\vbox to 5cm {}

\noindent avec pour $p=p_j$, $p^-:=p_{j-1}$ quand $1\le j\le k$, et $p^+:= \Big|\begin{array}{cl}
p_{j+1} &(0\le j<k)\\
p_\infty &(j=k).
\end{array}$ La notation $qa$ signifiera toujours implicitement que $q=P(aq)$, autrement dit que $q$ est un nombre premier et que
$$
P(a)\le q.
$$
Précisons que contrairement à ce que le schéma (8.16) suggère, certains des connecteurs $qa$ seront intérieurs et non extérieurs.

\vspace{2mm}

\noindent\textbf{Etape 2. Choix des connecteurs $qa$ pour les $p\in \{p_o,p_k\}$ ou $p\le 11$}

Cela représente au plus 6 nombres premiers $p$. On choisit les $qa$ correspondant comme des puissances de 2 deux à deux distinctes, non présentes dans $\Cc_{A,A'}\sqcup \Cc_{B,B'}$, et vérifiant
$$
\sqrt{x }< qa \le 2^8\sqrt{x}.
$$
C'est possible d'après (8.7) et (8.8).

Récapitulons en notant $Q\CA_2$ l'ensemble de ces connecteurs $qa$. On a
\begin{equation}
qa \in Q\CA_2 \Longrightarrow q=2
\end{equation}
\begin{equation}
qa\in Q\CA_2 \Longrightarrow qa \in ]\sqrt{x}, 2^8\sqrt{x}] \cap \CA^*(x,y)
\end{equation}
et
\begin{equation}
qa \in Q\CA_2 \Longrightarrow qa \notin \Cc_{A,A'}\cup \Cc_{B,B'}.
\end{equation}

Supposons qu'il existe des nombres premiers $p$ auxquels on n'a pas encore attribué de connecteur $qa$. C'est alors l'objet de l'étape suivante de les choisir.

\vspace{2mm}

\noindent\textbf{Etape 3. Choix des connecteurs $qa$ pour les $p$ vérifiant $13\le p\le y$ et $p\notin \{p_o,p_k,p_\infty\}$}

On note $p_{j_o}$ le plus petit de ces nombres premiers $p$. On choisit alors les connecteurs $qa=q_ja_j$ avec $j_o\le j<k$ fournis par le lemme~7.1 avec
$$
\Ee= \{p_o,p_\infty\} \cap [13,y].
$$
Notons $Q\CA_{13}$ l'ensemble de ces connecteurs. On a d'après le lemme~7.1
\begin{align}
qa\in Q\CA_{13} &\Longrightarrow q\ge 13\\
qa\in Q\CA_{13} &\Longrightarrow qa >\sqrt{x}\\
qa\in Q\CA_{13} &\Longrightarrow q\le p\\
qa\in Q\CA_{13} &\Longrightarrow 2^{13}p q^2\le x\ \mathrm{et\ }8pqS(a) \le x\\
qa\in Q\CA_{13} &\Longrightarrow q\notin \{p_o,p_\infty\}\\
qa\in Q\CA_{13} &\Longrightarrow a\le \sqrt{x/q}.
\end{align}

\eject

\noindent\textbf{Etape 4. Choix des extrémités des chaînes  $\Cc^*(x/p,p)$}

Rappelons que l'on a $p_o=P(A')$ et $p_\infty =P(B')$. Le premier entier de $\Cc_{A',B'}^*$ est $A'$. Donc l'extrémité gauche de $\Cc^*(x/p_o,p_o)$ est nécessairement $A'/p_o$.
De même, le dernier entier de $\Cc_{A',B'}^*$ est $B'$ et l'extrémité droite de $\Cc^*(x/p_\infty,p_\infty)$ est nécessairement $B'/p_\infty$.

On choisit à présent l'extrémité droite de $\Cc^*(x/p_o,p_o)$ comme étant $d_oa_oq_o$ avec $d_o\in \{1,3\}$ de telle sorte que
$$
P(d_o a_o q_o)\not= P(A'/p_o).
$$
De même on choisit l'extrémité gauche de $\Cc^*(x/p_\infty,p_\infty)$ comme étant $d_\infty a_k q_k$ avec $d_\infty \in \{1,3\}$ de telle sorte que
$$
P(d_\infty a_k q_k) \not= P(B'/p_\infty).
$$

Pour les autres $j$ (s'il en existe) qui vérifient $1\le j\le 4$ et $p=p_j\le \min(y,11)$, on choisit $q_{j-1}a_{j-1}$ comme extrémité gauche de $\Cc^*(x/p_j,p_j)$ et $3q_ja_j$ pour son extrémité droite.

S'il existe encore d'autres $p$ vérifiant $3\le p\le y$, on choisit $q(p^-)a(p^-)$ pour l'extrémité gauche de $\Cc^*(x/p,p)$ et $q(p)a(p)$ pour l'extrémité droite.

Tous les choix ont été faits. Il reste à procéder à plusieurs vérifications, dont les deux premières vont permettre de construire les chaînes $\Cc^*(x/p,p)$.

\vspace{2mm}

\noindent\textbf{Etape 5. Vérifications}

\vspace{2mm}

\textbf{Vérification 1.}

Il résulte de (8.12), (8.17) et (8.18) que
$$
\left|
\begin{array}{lrc}
&dqa \in \CA^*(x/p,p) \cap \CA^*(x/p^+,p^+)\\
\mathrm{et}\\
&dq_ka_k \in \CA^*(x/p_\infty,p_\infty)
\end{array}
\ \begin{pmatrix}
3 \le p\le y, \ p\not= p_\infty\\
d\in \{1,3\}\\
qa\in Q\CA_2
\end{pmatrix}
\right.
$$

Quand $p=p_j$ avec $1\le j<k$, il y a au plus deux nombres premiers (qui sont le cas échéant $p_o$ et/ou $p_\infty$) entre $p$ et $p^+$, d'où avec le postulat de Bertrand, $p^+<8p$. Avec  (8.21), (8.22) et (8.23), on a donc
$$
qa\in \CA^*(x/p,p) \cap \CA^*(x/p^+,p^+),\ \begin{pmatrix}
p=p_j\ \mathrm{avec\ }1\le j<k\\
qa\in Q\CA_{13}
\end{pmatrix}.
$$

\eject

\textbf{Vérification 2.}

Soit $p\in [3,y]$. On note $e_{p^-}$ et $e_p$ les deux extrémités de $\Cc^*(x/p,p)$ que l'on a choisies. Vérifions que
$$
P(e_{p^-}) \not= P(e_p), \qquad (3\le p\le y).
$$
On a fait en sorte que cela soit vrai à l'étape 4 quand $p\in \{p_o,p_\infty\}$ ou quand les deux connecteurs $qa$ qui sont voisins de $\Cc^*(x/p,p)$ appartiennent à $Q\CA_2$. Quand ces deux connecteurs appartiennent à $Q\CA_{13}$, cela résulte de la propriété (7.3) du lemme~7.1. Enfin quand ces deux connecteurs appartiennent l'un à $Q\CA_2$ et l'autre à $Q\CA_{13}$, cela résulte de (8.17) et (8.20).

En utilisant la définition (8.5) de $h_a^*(x,y)$  on déduit de ces deux premières vérifications que pour tout $p\in [3,y]$, on peut construire une chaîne $\Cc^*(x/p,p)$ d'entiers de $\CA^*(x/p,p)$, dont les extrémités sont  comme on les a choisies et qui vérifie
$$
\mathrm{longueur\ }(\Cc^*(x/p,p)) \ge h_a^*(x/p,p).
$$
En utilisant (8.14), on a donc pour $\Cc_{A',B'}^*$ et $\Cc_{A,B}^*$ les suites finies d'entiers que l'on a construites
$$
\mathrm{longueur}(\Cc^*_{A,B}) \ge \mathrm{longueur} (\Cc_{A',B'}^*) \ge 1+ \sum_{2<p\le y} h_a^*(x/p,p).
$$
Il reste à vérifier que $\Cc_{A,B}^*$ est une chaîne d'entiers de $\CA^*(x,y)$.

\vspace{2mm}

\textbf{Vérification 3.}

D'abord on a collé les chaînes $\Cc_{A,A'}^*$ et $\Cc_{B',B}^*$ à $\Cc_{A',B'}^*$. Ensuite on a juxtaposé les chaînes $p\Cc^*(x/p,p)$ par les connecteurs $qa$. Tout cela montre que pour toute paire de voisins de la suite $\Cc_{A,B}^*$, le petit divise le second. 

\textbf{Vérification 4.}

Soit $p\in [3,y]$, $p\not= p_\infty$.

D'une part avec (8.18) et (8.21), on a $qa>\sqrt{x}$.

Et d'autre part, par la vérification 1, on a $qa\in \CA^*(x/p,p)$.

Tout cela montre que
$$
qa\in \CA^*(x,y),\qquad (3\le p \le y,\ p\not= p_\infty).
$$

On a aussi
$$
p\Cc^*(x/p,p) \subset p \CA^*(x/p,p) \subset \CA^*(x,y)
$$
 d'après (8.12). Le lemme 8.6 permet alors de conclure que
$$
\Cc_{A,B}^* \subset \CA^*(x,y).
$$ 

Il reste à vérifier que les entiers de $\Cc_{A,B}^*$ sont deux à deux distincts. Répartissons les entiers de $\Cc_{A,B}^*$ en trois multiensembles
$$
\begin{array}{l}
E_1 = \Cc_{A,A'}^* \sqcup \Cc_{B',B}^*\\
E_2 = Q\CA_2 \sqcup Q\CA_{13}\\
E_3 = \bigsqcup\limits_{3\le p\le y} p\Cc^*(x/p,p).
\end{array}
$$
Il s'agit de vérifier d'une part que dans chacun de ces trois multiensembles d'entiers, les entiers sont deux à deux distincts, et donc que ce sont des ensembles, et d'autre part que 
$$
E_1 \cap E_2 = \emptyset,E_1\cap E_3= \{A',B'\},\  \mathrm{et}\ E_2\cap E_3=\emptyset.
$$
Faisons--le.

\vspace{2mm}

\textbf{Vérification 5.}

C'est le lemme 8.6 qui assure que les éléments de $E_1$ sont deux à deux distincts.

Les connecteurs choisis à l'étape 2 sont deux à deux distincts. Il découle du lemme 7.1 que les connecteurs choisis à l'étape 3 sont aussi deux à deux distincts. Enfin d'après (8.17) et (8.20), on a $Q\CA_2 \cap Q\CA_{13}=\emptyset$. Donc les éléments de $E_2$ sont deux à deux distincts. C'est (8.15) qui permet d'assurer que les éléments de $E_3$ sont deux à deux distincts. La propriété~(8.19) assure que $E_1\cap Q\CA_2=\emptyset$. La combinaison de la propriété~(8.10) du lemme~8.6 avec les propriétés (8.20), (8.24), (8.22) et (8.12)  montrent que $E_1\cap Q\CA_{13}=\emptyset$. Donc finalement,
$$
E_1\cap E_2=\emptyset.
$$

Soient $p\in [3,y]$ et $m\in \Cc^*(x/p,p)$. On a  $P(pm)=p$. En utilisant (8.12) et l'inclusion $\Cc^*(x/p,p) \subset \CA^*(x/p,p)$, on en déduit que
$$
\begin{array}{c}
3 \le P(pm) \le \sqrt{x} /2^{\gamma+13},\\
\vspace{2mm}
pm >\sqrt{xP(pm)}
\end{array}
$$
et
$$
P(pm/P(pm)) \le \sqrt{x/P(pm)}/27.
$$
Avec (8.11), on en déduit que
$$
E_1 \cap E_3 = \{A',B'\}.
$$

D'après (8.17), on a $Q\CA_2\cap E_3 =\emptyset$. D'après (8.25), on a
\begin{equation}
qa \in Q\CA_{13} \Longrightarrow qa \le \sqrt{xP(qa)}.
\end{equation}
Or on a par ailleurs
\begin{equation}
m \in \Cc^*(x/p,p) \Longrightarrow pm > \sqrt{xP(pm)},
\end{equation}
d'où $Q\CA_{13}\cap E_3 = \emptyset$, soit finalement
$$
E_2 \cap E_3 = \emptyset.
$$ 

On a achevé les vérifications. Cela conclut  la preuve du lemme~8.7.

\vspace{2mm}

\noindent\textbf{Remarque.} C'est la combinaison de (8.26) et (8.27) qui nous a permis de vérifier facilement que les entiers de l'étage du bas de (8.16) sont différents de ceux du haut. C'est  cette vérification facile qui nous a conduit, pour établir la minoration du théorème~1.3, d'imposer aux entiers $n$ de la chaîne d'entiers de $\CA(x,y)$ que l'on construit, de vérifier la condition supplémentaire $n>\sqrt{x}$.

\section{Etude de la fonction $h_a^*(x,y)$}

{\monlem Pour tout $x\ge 3^{10}$ fixé, l'application
$$
\begin{array}{cll}
[3,+\infty[ &\longrightarrow &\N^*\\
y & \longmapsto &h_a(x,y)
\end{array}
$$
est croissante}

\vspace{2mm}

\noindent\textbf{Démonstration.} Supposons $3\le y <z$.

Soit $(A,B)\in \CA^*(x,z)^2$ avec $P(A)<P(B)$. Si $(A,B) \in \CA^*(x,y)^2$, alors on a immédiatement
$$
f_{a,A,B}^* (x,z) \ge f_{a,A,B}^*(x,y)
$$

Si $A\notin \CA^*(x,y)$, alors $y<P(A)<P(B)\le z$. En appliquant le lemme~8.1, on obtient une chaîne $\Cc_{A,2^\alpha}^*$ d'entiers $A=a_1-a_2-\cdots - a_{a_s}=2^\alpha$ avec $P(a_j)=P(A)$ pour tout $j$ vérifiant $1\le j<s$.

En appliquant le lemme 8.2, on obtient une chaîne $\Cc_{3^\beta,B}$ d'entiers $3^{\beta} = b_1 - b_2 -\cdots - b_{s'}=B$ avec $P(b_j)=P(B)$ pour tout $j$ vérifiant $2\le j\le s'$. On choisit alors une chaîne $\Cc_{2^\alpha,3^{\beta}}^*$ d'entiers de $\CA^*(x,y)$ et de longueur $f_{a,2^\alpha,3^{\beta}}^*(x,y)$. Alors
$$
\Cc_{A,2^\alpha}^* - \Cc_{2^\alpha,3^{\beta}}^* - \Cc_{3^{\beta},B}^*
$$
forme une chaîne d'entiers de $\CA^*(x,z)$ de longueur $>f_{a,2^\alpha,3^{\beta}}(x,y)$. D'où
$$
f_{a,A,B}^*(x,z) > f_{a,2^\alpha,3^{\beta}}^*(x,y).
$$
Si enfin $P(A) \le y < P(B)\le z$, c'est avec une chaîne $\Cc_{A,2^{\alpha}}^* - \Cc_{2^{\alpha},B}^*$ que l'on montre que $f_{a,A,B}(x,z) > f_{a,A,2^{\alpha}}^*(x,y)$. Tout cela montre que $h(x,z) \ge h(x,y)$.

\vspace{2mm}

{\monlem Il existe un réel positif $x_2$ et un entier positif $\gamma$ tels que pour tous réels $x$ et $y$ vérifiant
$$
x\ge x_2 \qquad \mathrm{et}\ y\ge 3
$$
on a
$$
h_a^*(x,y) \ge 1 + \sum_{2<p\le\min(y,\sqrt{x}/2^{\gamma +13})} h_a^*(x/p,p).
$$
}

\vspace{2mm}

\noindent\textbf{Démonstration.} En combinant les lemmes 8.5 et 8.7, on sait qu'il existe un entier positif $\gamma$ et un réel positif $x_2$ tels que pour $3\le y\le \sqrt{x}/2^{\gamma+13}$ et $x\ge x_2$, on a
$$
h_a^* (x,y) \ge 1 + \sum_{2<p\le y}h_a^*(x/p,p).
$$
En utilisant le lemme 9.1, on en déduit que pour $x\ge x_0$ et $y\ge 3$, on a
$$
\begin{array}{rl}
h_a^*(x,y) &\ge h_a^*(x,\min(y,\sqrt{x} /2^{\gamma+13}))\\
&\ge 1+ \displaystyle \sum_{2<p\le\min(y,\sqrt{x}/2^{\gamma+13})} h_a^*(x/p,p).
\end{array}
$$
Cela achève la preuve du lemme 9.2.

\section{Preuve du théorème 1.4}

Soient $x$ et $y$ des réels vérifiant $x>0$ et $y\ge 2$. Posons
$$
M=\max (2^{2(\gamma+13)},x_2/9)
$$
où $\gamma$ et $x_2$ sont convenables pour le lemme 9.2. D'après le lemme 9.2, on a
\setcounter{equation}{0}
\begin{equation}
h_a^*(x,y) \ge 1\!\!1_{[1,+\infty[} + \sum_{2<p\le\min(y,\sqrt{x/M})}h_a^*(x/p,p)
\end{equation}
sous les conditions
$$
x\ge x_2\quad \mathrm{et}\quad y\ge 3.
$$
Si ces conditions ne sont pas vérifiées, alors la somme en $p$ dans (10.1) est vide et l'inégalité (10.1) est immédiate. Elle est donc satisfaite dès que $x>0$ et $y\ge 2$. En appliquant le lemme~5.1, on a donc pour tous réels $x$ et $y$ vérifiant $x\ge 8$ et $y\ge 2$,
$$
f_a(x,y) \ge f_a^*(x,y) \ge h_a^*(x,y) \ge A(x/M,y,2).
$$
Cela achève la preuve du théorème 1.4.

\eject

\eject

\vskip4mm

\begin{tabular}{ll}

 &\hspace{6.7cm}Eric Saias \\

 &\hspace{6.7cm}Sorbonne Université\\

&\hspace{6.7cm}LPSM\\

&\hspace{6.7cm}4, place Jussieu\\

&\hspace{6.7cm}75252 Paris Cedex 05 (France)\\

\vspace{2mm}

 &\hspace{6.7cm}\textsf{eric.saias@upmc.fr}
\end{tabular}

\end{document}